\documentclass[12pt]{article}
\usepackage{amsmath}
\usepackage{amssymb}
\usepackage{vatola}

\def\q{\quad}
\def\qq{\qquad}

\def\mod#1{\ (\text{\rm mod}\ #1)}
\def\t{\text}
\def\f{\frac}
\def\e{\equiv}
\def\b{\binom}

\def\sls#1#2{(\f{#1}{#2})}
 \def\ls#1#2{\big(\f{#1}{#2}\big)}
\def\Ls#1#2{\Big(\f{#1}{#2}\Big)}

\let \pro=\proclaim
\let \endpro=\endproclaim

\begin{document}
 \centerline {\bf
New conjectures involving binomial coefficients and Ap\'ery-like
numbers}
\par\q\newline
\centerline{Zhi-Hong Sun}\newline \centerline{School of Mathematics
and Statistics}\centerline{Huaiyin Normal University}
\centerline{Huaian, Jiangsu 223300, P.R. China} \centerline{Email:
zhsun@hytc.edu.cn} \centerline{Homepage:
http://maths.hytc.edu.cn/szh1.htm}
 \abstract{In this paper, we pose lots of challenging conjectures on
 congruences for the sums involving binomial coefficients
  and Ap\'ery-like numbers modulo $p^3$, where $p$ is an odd prime.
 \par\q
 }
 \endabstract
\section*{1. Introduction}
\par\q For $a\in\Bbb Z$ and given
odd prime $p$ let $\sls ap$ denote the Legendre symbol.  For
positive integers $a,b$ and $n$, if $n=ax^2+by^2$ for some integers
$x$ and $y$, we briefly write that $n=ax^2+by^2$. Let $p>3$ be a
prime. In 1987, Beukers[B] conjectured a congruence equivalent to
$$\sum_{k=0}^{p-1}\f{\b{2k}k^3}{64^k}\e
\cases 4x^2-2p\mod{p^2}&\t{if $p=x^2+4y^2\e 1\mod 4$,}
\\ 0\mod{p^2}&\t{if $p\e 3\mod 4$.}\endcases$$
This congruence was proved by several authors including
Ishikawa[I]($p\e 1\mod 4$), van Hamme[vH]($p\e 3\mod 4$) and
Ahlgren[A]. Combining the results in [LR], [S6] and [T], in [S11]
the author stated that
$$\sum_{k=0}^{p-1}\f{\b{2k}k^3}{64^k}\e
\cases 4x^2-2p-\f{p^2}{4x^2}\mod{p^3}&\t{if $p=x^2+4y^2\e 1\mod
4$,}\\-\f {p^2}4\b{(p-3)/2}{(p-3)/4}^{-2}\mod{p^3}&\t{if $p\e 3\mod
4$.}\\\endcases$$

\par Let $p$ be a prime with $p\not=2,7$. In 1998,  using the hypergeometric series $\,_3F_2(\lambda)_p$
over the finite field $\Bbb F_p$, Ono[O] obtained some congruences
equivalent to $\sum_{k=0}^{p-1}\f 1{m^k}\b{2k}k^3$ modulo $p$ in the
cases $m=1,-8,16,$ $-64,256,-512,4096$. For such values of $m$, in
[Su1,Su2] the author's brother Zhi-Wei Sun conjectured the
congruences for $\sum_{k=0}^{p-1}\f 1{m^k}\b{2k}k^3$ modulo $p^2$,
which have been proved in [S2] and [KLMSY]. In [S7], the author
conjectured the congruences for $\sum_{k=0}^{p-1}\f 1{m^k}\b{2k}k^3$
modulo $p^3$. For example, for any prime $p\not=2,3,7$,
$$\sum_{k=0}^{p-1}\b{2k}k^3\e
\cases  4x^2-2p-\f{p^2}{4x^2}\mod {p^3}\q \t{if $p\e 1,2,4\mod 7$
and so $p=x^2+7y^2$,}
\\-\f {11}4p^2\b{3[p/7]}{[p/7]}^{-2} \e
-11p^2\b{[3p/7]}{[p/7]}^{-2} \mod {p^3}\qq\ \t{if $7\mid p-3$,}
\\-\f{99}{64}p^2\b{3[p/7]}{[p/7]}^{-2}\e
-\f{11}{16}p^2\b{[3p/7]}{[p/7]}^{-2}\mod {p^3}\qq\ \t{if $7\mid
p-5$,}
\\-\f{25}{176}p^2\b{3[p/7]}{[p/7]}^{-2}\e
-\f{11}4p^2\b{[3p/7]}{[p/7]}^{-2}\mod {p^3}\qq\t{if $7\mid p-6$,}
\endcases$$
where $[x]$ is the greatest integer not exceeding $x$. In [S11], the
author posed conjectures on
$$\sum_{k=0}^{p-1}
\f{\b{2k}k^3}{m^k(k+1)} \q \sum_{k=0}^{p-1}
\f{\b{2k}k^3}{m^k(2k-1)}, \q \sum_{k=0}^{p-1}
\f{\b{2k}k^3}{m^k(2k-1)^2}, \q\sum_{k=0}^{p-1}
\f{\b{2k}k^3}{m^k(2k-1)^3}$$
 modulo $p^3$, and $\sum_{k=0}^{p-1}
\f{\b{2k}k^3}{m^k(k+1)^2}$ modulo $p^2$. For instance, for any odd
prime $p\e 1,2,4\mod 7$ and so $p=x^2+7y^2$,
$$\align &\sum_{k=0}^{(p-1)/2}
\f{\b{2k}k^3}{k+1}\e -44y^2+2p\mod {p^3},
\\&\sum_{k=0}^{p-1}
\f{\b{2k}k^3}{(k+1)^2}\e -68y^2\mod {p^2},
\\&\sum_{k=0}^{p-1}
\f{\b{2k}k^3}{2k-1}\e -36y^2+14p-\f{7p^2}{4y^2}\mod {p^3},
\\&\sum_{k=0}^{p-1}
\f{\b{2k}k^3}{(2k-1)^2}\e -284y^2+34p+\f{23p^2}{4y^2}\mod {p^3},
\\&\sum_{k=0}^{p-1}
\f{\b{2k}k^3}{(2k-1)^3}\e -804y^2-18p-\f{39p^2}{4y^2}\mod
{p^3}.\endalign$$

Let $p>3$ be a prime. In 2003, Rodriguez-Villegas[RV] posed 22
conjectures on supercongruences modulo $p^2$. In particular, the
following congruences are equivalent to conjectures due to
Rodriguez-Villegas:
$$\align &\sum_{k=0}^{p-1}\f{\b{2k}k^2\b{3k}k}{108^k}\e
\cases 4x^2-2p\mod{p^2}&\t{if $p=x^2+3y^2\e 1\mod 3$,}
\\0\mod{p^2}&\t{if $p\e 2\mod 3$,}
\endcases
\\&\sum_{k=0}^{p-1}\f{\b{2k}k^2\b{4k}{2k}}{256^k}\e
\cases 4x^2-2p\mod{p^2}&\t{if $p=x^2+2y^2\e 1,3\mod 8$,}
\\0\mod{p^2}&\t{if $p\e 5,7\mod 8$,}
\endcases
\\&\Ls p3\sum_{k=0}^{p-1}\f{\b{2k}k\b{3k}k\b{6k}{3k}}{12^{3k}}\e
\cases 4x^2-2p\mod{p^2}&\t{if $p=x^2+4y^2\e 1\mod 4$,}
\\0\mod{p^2}&\t{if $p\e 3\mod 4$.}
\endcases\endalign$$
  These conjectures have been solved  by
Mortenson[M1] and Zhi-Wei Sun[Su3]. In 2018, J.C. Liu[Liu1]
conjectured the congruences for
$$\sum_{k=0}^{p-1}\f{\b{2k}k^2\b{3k}k}{108^k},\q
\sum_{k=0}^{p-1}\f{\b{2k}k^2\b{4k}{2k}}{256^k},\q
\sum_{k=0}^{p-1}\f{\b{2k}k\b{3k}k\b{6k}{3k}}{1728^k}\mod {p^3}$$ in
terms of $p-$adic Gamma functions.

Let $p$ be an odd prime, $m\in\Bbb Z$ and $p\nmid m$. In [Su1,Su4],
  Z.W. Sun posed many conjectures concerning congruences modulo $p^2$
involving the sums
$$\sum_{k=0}^{p-1}\f{\b{2k}k^2\b{3k}k}{m^k},\q
\sum_{k=0}^{p-1}\f{\b{2k}k^2\b{4k}{2k}}{m^k}, \q
\sum_{k=0}^{p-1}\f{\b{2k}k\b{3k}k\b{6k}{3k}}{m^k}.$$ For 13 similar
conjectures see [S1].  Most of these congruences modulo $p$ were
proved by the author in [S2-S5]. In [S7] and [S11], the author
conjectured many congruences for
$$\align &\sum_{k=0}^{p-1}\f{\b{2k}k^2\b{3k}k}{m^k},\q
\sum_{k=0}^{p-1}\f{\b{2k}k^2\b{3k}k}{m^k(k+1)},\q
\sum_{k=0}^{p-1}\f{\b{2k}k^2\b{3k}k}{m^k(2k-1)},\q
\\&\sum_{k=0}^{p-1}\f{\b{2k}k^2\b{4k}{2k}}{m^k}, \q
\sum_{k=0}^{p-1}\f{\b{2k}k^2\b{4k}{2k}}{m^k(k+1)}, \q
\sum_{k=0}^{p-1}\f{\b{2k}k^2\b{4k}{2k}}{m^k(2k-1)}, \q
\\&\sum_{k=0}^{p-1}\f{\b{2k}k\b{3k}k\b{6k}{3k}}{m^k},\q
\sum_{k=0}^{p-1}\f{\b{2k}k\b{3k}k\b{6k}{3k}}{m^k(k+1)},\q
\sum_{k=0}^{p-1}\f{\b{2k}k\b{3k}k\b{6k}{3k}}{m^k(2k-1)}\endalign$$
modulo $p^3$. As typical examples, for any prime $p>5$,
$$\align &\sum_{k=0}^{p-1}\f{\b{2k}k^2\b{3k}k}{(-192)^k(k+1)}
\e\cases \f 32x^2-4p-p^2\mod{p^3}\q\t{if $3\mid p-1$ and so
$4p=x^2+27y^2$,}
\\2(2p+1)\b{[2p/3]}{[p/3]}^2+p\mod {p^2}\q\t{if $p\e 2\mod 3$,}
\endcases
\\&\sum_{k=0}^{p-1}\f{\b{2k}k^2\b{3k}{k}}{(-192)^k(2k-1)} \\&\q\e\cases -\f
34x^2+\f {9p}8+\f{3p^2}{8x^2}\mod{p^3}&\t{if $3\mid p-1$ and so
$4p=x^2+27y^2$,}
\\-\f 12(2p+1)\b{[2p/3]}{[p/3]}^2+\f 38p\mod {p^2}&\t{if $p\e 2\mod{3}$,}
\endcases
\\&\sum_{k=0}^{p-1}\f{\b{2k}k^2\b{4k}{2k}}{648^k(k+1)} \e\cases -\f
{40}3y^2+2p-p^2\mod {p^3}&\t{if $p=x^2+4y^2\e 1\mod 4$,}
\\-\f 32R_1(p)-\f p3\mod {p^2}&\t{if $p\e 3\mod 4$,}
\endcases
\\&\sum_{k=0}^{p-1}\f{\b{2k}k^2\b{4k}{2k}}{648^k(2k-1)}
\e \cases -\f{76}{27}x^2+\f{104}{81}p+\f{67}{324x^2}p^2\mod
{p^3}&\t{if $p=x^2+4y^2\e 1\mod 4$,}
\\-\f 29R_1(p)+\f{10}{81}p
\mod {p^2}&\t{if $p\e 3\mod 4$,}
\endcases
\\&\Ls p5 \sum_{k=0}^{p-1}\f{\b{2k}k\b{3k}{k}\b{6k}{3k}}{54000^k(k+1)}
\e\cases \f{48}5y^2+2p-\ls p5p^2\mod {p^3} &\t{if $p=x^2+3y^2\e
1\mod 3$,}
\\-20R_3(p)-\f{18}5p\mod {p^2}&\t{if $p\e 2\mod 3$,}
\endcases
\\&\Ls p5\sum_{k=0}^{p-1}\f{\b{2k}k\b{3k}{k}\b{6k}{3k}}{54000^k(2k-1)}
\\&\q\e\cases -\f {748}{225}x^2+\f{1708}{1125}p+\f{103}{500x^2}p^2\mod
{p^3}&\t{if $p=x^2+3y^2\e 1\mod 3$,}
\\-\f 8{45}R_3(p)+\f{18}{125}p\mod {p^2}&\t{if $p\e 2\mod 3$,}
\endcases
\endalign$$
where
$$\align &R_1(p)=(2p+2-2^{p-1})\b{(p-1)/2}{[p/4]}^2,
\\&R_3(p)=\Big(1+2p+\f
43(2^{p-1}-1)-\f 32(3^{p-1}-1)\Big) \b{(p-1)/2}{[p/6]}^2.
\endalign$$

\par
 In Section 2, with the help of Maple, we pose new
conjectures on congruences modulo $p^3$ involving the sums
$$\align
&\sum_{k=0}^{p-1}k^2\f{\b{2k}k^3}{m^k},\q
\sum_{k=0}^{p-1}k^3\f{\b{2k}k^3}{m^k},\q
\sum_{k=0}^{p-1}k^2\f{\b{2k}k^2\b{3k}k}{m^k},\q
\sum_{k=0}^{p-1}k^3\f{\b{2k}k^2\b{3k}k}{m^k},
\\&\sum_{k=0}^{p-1}k^2\f{\b{2k}k^2\b{4k}{2k}}{m^k}, \q
\sum_{k=0}^{p-1}k^3\f{\b{2k}k^2\b{4k}{2k}}{m^k},\q
\sum_{k=0}^{p-1}k^2\f{\b{2k}k\b{3k}k\b{6k}{3k}}{m^k},\q
\sum_{k=0}^{p-1}k^3\f{\b{2k}k\b{3k}k\b{6k}{3k}}{m^k}.
\endalign$$
For instance, for any odd prime $p=x^2+7y^2\e 1,2,4\mod 7$,
$$\align&\sum_{k=0}^{p-1}k^2\b{2k}k^3\e
 \f{736x^2}{1323}-\f{272p}{441}+\f{20p^2}{1323x^2}\mod {p^3},
\\&\sum_{k=0}^{p-1}k^3\b{2k}k^3\e
-\f{5408x^2}{27783}+\f{2992p}{9261}-\f{1774p^2}{27783x^2}\mod {p^3}.
\endalign$$
We remark that Z.W. Sun[Su1] proved for any odd prime $p$,
$$\sum_{k=0}^{p-1}(21k+8)\b{2k}k^3\e 8p\mod {p^3}.$$
\par The first kind of Ap\'ery-like
numbers $\{u_n\}$ satisfies
$$u_0=1, \ u_1=\ b,\ (n+1)^3u_{n+1}=(2n+1)(an(n+1)+b)u_n-cn^3u_{n-1}\q (n\ge 1),
$$ where $a,b,c\in\Bbb Z$, $c\not=0$ and $u_n\in\Bbb Z$ for all
positive integers $n$. Let
$$\align
 &A_n= \sum_{k=0}^n\binom nk^2\binom{n+k}k^2,
\\& D_n=\sum_{k=0}^n\b nk^2\b{2k}k\b{2n-2k}{n-k},\\&
b_n=\sum_{k=0}^{[n/3]}\b{2k}k\b{3k}k\b n{3k}\b{n+k}k(-3)^{n-3k}, \\&
T_n=\sum_{k=0}^n\b
nk^2\b{2k}n^2=\sum_{k=0}^{[n/2]}\b{2k}k^2\b{4k}{2k}\b{n+2k}{4k}4^{n-2k},
\
\\&V_n=\sum_{k=0}^n
\b
nk\b{n+k}k(-1)^k\b{2k}k^216^{n-k}=\sum_{k=0}^n\b{2k}k^2\b{2n-2k}{n-k}^2
\\& \q=\sum_{k=0}^n\b{2k}k^3\b k{n-k}(-16)^{n-k},
\\&V_n^{(3)}=\sum_{k=0}^n\b
nk\b{n+k}k(-1)^k\b{2k}k\b{3k}k27^{n-k},
\\&V_n^{(4)}=\sum_{k=0}^n\b
nk\b{n+k}k(-1)^k\b{2k}k\b{4k}{2k}64^{n-k}
=\sum_{k=0}^n\b{2k}k^3\b{2n-2k}{n-k}16^{n-k},
\\&V_n^{(6)}=\sum_{k=0}^n\b
nk\b{n+k}k(-1)^k\b{3k}k\b{6k}{3k}432^{n-k}.\endalign$$
 Then $\{A_n\}$, $\{D_n\}$, $\{b_n\}$, $\{T_n\}$, $\{V_n\}$, $\{V_n^{(3)}\}$,
  $\{V_n^{(4)}\}$ and $\{V_n^{(6)}\}$ are
 the first kind of Ap\'ery-like numbers with
$(a,b,c)=(17,5,1),\ (10,4,64),\ (-7,-3,81),\ (12,4,16)$,
$(16,8,256)$,\newline $(27,15,729),$ $(64,40,4096)$ and
$(432,312,186624)$, respectively. The numbers $\{A_n\}$, $\{D_n\}$
and $\{b_n\}$ are called Ap\'ery numbers, Domb numbers and
Almkvist-Zudilin numbers, respectively. For $\{A_n\}$, $\{D_n\}$,
$\{b_n\}$, $\{T_n\}$ and $\{V_n\}$ see A005259, A002895, A125143,
A290575 and A036917 in Sloane's database ``The On-Line Encyclopedia
of Integer Sequences". For the congruences concerning $\{A_n\}$,
$\{D_n\}$ and $\{b_n\}$ see [Su2,Su4,S5]. For the congruences
involving $T_n$ see the author's paper [S6]. For the formulas and
congruences involving $V_n$ see [AZ,S8,S9,Su5,W,Z]. For the
congruences concerning $V_n^{(3)},\;V_n^{(4)}$ and $V_n^{(6)}$ see
[S9]. In [S6-S11], the author conjectured many congruences modulo
$p^3$ involving Ap\'ery-like numbers.
 \par Let $p>3$ be a prime, $m\in\Bbb Z$ and $p\nmid m$. In Section
3, we pose many conjectures on $\sum_{k=0}^{p-1}\f{k^2u_k}{m^k}$ and
$\sum_{k=0}^{p-1}\f{k^3u_k}{m^k}$ modulo $p^3$, where
$u_n\in\{A_n,D_n,b_n,T_n,V_n,V_n^{(3)},V_n^{(4)},V_n^{(6)}\}$. For
example,
$$
\align &\sum_{n=0}^{p-1}n^2\f{T_n}{4^n} \e\cases -4y^2\mod
{p^3}&\t{if $p=x^2+4y^2\e 1\mod 4$,}
\\-\f 14R_1(p)-\f 12p\mod {p^2}
&\t{if $p\e 3\mod 4$,}
\endcases
\\&\sum_{n=0}^{p-1}n^2\f{V_n^{(3)}}{243^n}
\e\cases \f {1}{16}x^2-\f {1}{8}p-\f{p^2}{8x^2}\mod{p^3}&\t{if
$3\mid p-1$ and so $4p=x^2+27y^2$,}
\\\f 18(2p+1)\b{[2p/3]}{[p/3]}^2\mod {p^2}&\t{if $p\e 2\mod{3}$,}
\endcases
\\&\sum_{n=0}^{p-1}n^3A_n \e\cases
-\f{13}{32}x^2+\f{37}{64}p-\f{19p^2}{256x^2}\mod {p^3}&\t{if
$p=x^2+2y^2\e 1,3\mod 8$,}
\\ \f 9{128}R_2(p)+\f 38p\mod {p^2}&\t{if $p\e 5,7\mod 8$,}\endcases
\endalign$$ where $R_2(p)$ is given by
$$\aligned R_2(p)=&(5-4(-1)^{\f{p-1}2})\Big(1+(4+2(-1)^{\f{p-1}2})p
-4(2^{p-1}-1)-\f p2\sum_{k=1}^{[p/8]} \f
1k\Big)\\&\q\times\b{\f{p-1}2}{[\f p8]}^2.\endaligned$$
\par For later convenience, we introduce the definitions of $\{B_n\},\{E_n\}$
and $\{U_n\}$. The Bernoulli numbers $\{B_n\}$, Euler numbers
$\{E_n\}$ and the sequence $\{U_n\}$ are defined by
$$\align &B_0=1,\q\sum_{k=0}^{n-1}\b nkB_k=0\q(n\ge 2),
\\& E_0=1,\q E_n=-\sum_{k=1}^{[n/2]}\b n{2k}E_{n-2k}\q(n\ge 1),
\\& U_0=1,\q U_n=-2\sum_{k=1}^{[n/2]}\b n{2k}U_{n-2k}\q(n\ge 1).
\endalign$$
It is known that $B_{2n+1}=E_{2n-1}=U_{2n-1}=0$ for $n\ge 1$.
\section*{2. Conjectures on congruences involving
binomial coefficients}

\par\q Calculations with Maple
suggest the following challenging conjectures:
\pro{Conjecture 2.1}
Let $p>3$ be a prime. Then
$$\align
&\sum_{k=0}^{p-1}\f{k^2\b{2k}k^2\b{3k}{k}}{(-192)^k}
\\&\e\cases \f
{2}{125}x^2-\f {27p}{250}+\f{11p^2}{250x^2}\mod{p^3}&\t{if $3\mid
p-1$ and so $4p=x^2+27y^2$,}
\\\f 2{25}(2p+1)\b{[2p/3]}{[p/3]}^2+\f {19}{250}p\mod {p^2}&\t{if $p\e 2\mod{3}$
and $p\not=5$,}
\endcases
\\&\sum_{k=0}^{p-1}\f{k^3\b{2k}k^2\b{3k}{k}}{(-192)^k}
\\&\e\cases \f
{21}{6250}x^2-\f {221p}{12500}-\f{447p^2}{12500x^2}\mod{p^3}&\t{if
$3\mid p-1$ and so $4p=x^2+27y^2$,}
\\-\f {27}{625}(2p+1)\b{[2p/3]}{[p/3]}^2+\f {137}{12500}p\mod {p^2}
&\t{if $p\e 2\mod{3}$ and $p\not=5$.}
\endcases
\endalign$$
\endpro

\pro{Conjecture 2.2} Let $p>5$ be a prime. Then
 $$\align
&\Ls{10}p\sum_{k=0}^{p-1}\f{k^2\b{2k}k\b{3k}{k}\b{6k}{3k}}
{(-12288000)^k}
\\&\e\cases \f 1{32388554}(\f {121213}2x^2-242919p+\f{493}{x^2}p^2)\mod{p^3}
\ \t{if $3\mid p-1$ and $4p=x^2+27y^2$,}
\\\f {800}{64009}(2p+1)\b{[2p/3]}{[p/3]}^2+\f{60853}{16194277}p
\mod {p^2}\q\ \t{if $3\mid p-2$ and $p\not=11,23$.}
\endcases
\endalign$$
\endpro
\par{\bf Remark 2.1} Let $p$ be a prime with $p>5$. In [S7], the
author conjectured that if $p\e 1\mod 3$ and so $4p=x^2+27y^2$, then
$$\sum_{k=0}^{p-1}\f{\b{2k}k^2\b{3k}k}{(-192)^k}\e
\Ls{10}p\sum_{k=0}^{p-1}\f{\b{2k}k\b{3k}k\b{6k}{3k}}{(-12288000)^k}
\e x^2-2p-\f{p^2}{x^2}\mod {p^3};$$ if $p\e 2\mod 3$, then
$$\sum_{k=0}^{p-1}\f{\b{2k}k^2\b{3k}k}{(-192)^k}
\e\f{800}{161}\Ls{10}p\sum_{k=0}^{p-1}\f{\b{2k}k\b{3k}k\b{6k}{3k}}{(-12288000)^k}
\e\f 34p^2\b{[2p/3]}{[p/3]}^{-2}\mod {p^3}.$$ The congruence for
$\sum_{k=0}^{p-1}\f{\b{2k}k^2\b{3k}k}{(-192)^k}$ modulo $p^2$ was
conjectured by Z.W. Sun[Su1] earlier, and solved in [S3] and [WS].
In [Su1], Z.W. Sun conjectured that
$$\sum_{k=0}^{p-1}(5k+1)\f{\b{2k}k^2\b{3k}k}{(-192)^k}\e\Ls p3p\mod
{ p^3}.$$ In [S5], the author conjectured that
$$\sum_{k=0}^{p-1}(506k+31)\f{\b{2k}k\b{3k}k\b{6k}{3k}}{(-12288000)^k}
\e 31\Ls{-30}pp\mod {p^3}.$$

\par Let $p>3$ be a prime. In 2008, Mortenson[M2] proved the
following congruence conjectured by van Hamme:
$$\sum_{k=0}^{p-1}(4k+1)\f{\b{2k}k^3}{(-64)^k}\e
(-1)^{\f{p-1}2}p\mod {p^3}.$$ In [GZ], Guillera and W. Zudilin
proved that
$$\sum_{k=0}^{p-1}(3k+1)\f{\b{2k}k^3}{(-8)^k}\e
(-1)^{\f{p-1}2}p\mod {p^3}.$$ In [Su1], Z.W. Sun conjectured that
$$\align &\sum_{k=0}^{p-1}(3k+1)\f{\b{2k}k^3}{(-8)^k}\e
(-1)^{\f{p-1}2}p+p^3E_{p-3}\mod {p^4},
\\&\sum_{k=0}^{(p-1)/2}(6k+1)\f{\b{2k}k^3}{(-512)^k}\e
\Ls {-2}pp+\f 14\Ls 2pp^3E_{p-3}\mod {p^4},
\\&\sum_{k=0}^{p-1}(3k+1)\f{\b{2k}k^3}{16^k}\e
p+\f 76p^4B_{p-3}\mod {p^5},
\\&\sum_{k=0}^{p-1}(6k+1)\f{\b{2k}k^3}{256^k}\e
(-1)^{\f{p-1}2}p-p^3E_{p-3}\mod {p^4},
\\&\sum_{k=0}^{p-1}(42k+5)\f{\b{2k}k^3}{4096^k}\e
5(-1)^{\f{p-1}2}p-p^3E_{p-3}\mod {p^4}.
\endalign$$

 \pro{Conjecture 2.3} Let $p$ be an odd prime. Then
$$\align
&\sum_{k=0}^{p-1}\f{k^2\b{2k}k^3}{(-8)^k} \e \cases \f{10}{27}x^2-\f
49p+\f{p^2}{54x^2}\mod {p^3}&\t{if $p=x^2+4y^2\e 1\mod 4$,}
\\\f 1{18}R_1(p)+\f 7{27}p
\mod {p^2}&\t{if $p\e 3\mod 4$ and $p\not=3$,}
\endcases
\\&\sum_{k=0}^{p-1}\f{k^3\b{2k}k^3}{(-8)^k} \e \cases -\f{4}{81}x^2+\f
4{27}p-\f{17p^2}{324x^2}\mod {p^3}&\t{if $p=x^2+4y^2\e 1\mod 4$,}
\\-\f 2{27}R_1(p)-\f {10}{81}p
\mod {p^2}&\t{if $p\e 3\mod 4$ and $p\not=3$.}
\endcases
\endalign$$
\endpro

\pro{Conjecture 2.4} Let $p$ be a prime of the form $4k+1$ and so
$p=x^2+4y^2$. Then
$$\sum_{k=0}^{p-1}\f{k^2\b{2k}k^3}{64^k}
\e  \f {x^2}{6}-\f p{12}-\f{p^2}{24x^2}\mod {p^3}.$$
\endpro
\par{\bf Remark 2.2} Let $p>3$ be a prime. In [S2], the author
obtained the congruences for $\sum_{k=0}^{p-1}\f{k\b{2k}k^3}{64^k},
\ \sum_{k=0}^{p-1}\f{k^2\b{2k}k^3}{64^k}$ and
$\sum_{k=0}^{p-1}\f{k^3\b{2k}k^3}{64^k}$ modulo $p^2$. See also [T].
\vskip0.2cm
 \pro{Conjecture 2.5} Let $p$ be an odd prime. Then
$$\align
&(-1)^{[\f p4]}\sum_{k=0}^{p-1}\f{k^2\b{2k}k^3}{(-512)^k}
\\&\q\e \cases \f {x^2}{27}-\f p{18}+\f {p^2}{216x^2}\mod {p^3}&\t{if
$p=x^2+4y^2\e 1\mod 4$,}
\\-\f 1{18}R_1(p)-\f p{27}\mod {p^2}
&\t{if $p\e 3\mod 4$ and $p\not=3$,}\endcases
\\&(-1)^{[\f p4]}\sum_{k=0}^{p-1}\f{k^3\b{2k}k^3}{(-512)^k}
\\&\q\e \cases -\f {x^2}{162}-\f p{216}-\f {p^2}{1296x^2}\mod {p^3}&\t{if
$p=x^2+4y^2\e 1\mod 4$,}
\\\f 1{108}R_1(p)-\f 5{648}p\mod {p^2}
&\t{if $p\e 3\mod 4$ and $p\not=3$.}
\endcases
\endalign$$
\endpro

\pro{Conjecture 2.6} Let $p>3$ be a prime. Then
$$\align
&\sum_{k=0}^{p-1}\f{k^2\b{2k}k^2\b{4k}{2k}}{648^k} \e \cases
\f{34}{343}x^2-\f{8p}{343}-\f{13p^2}{686x^2}\mod {p^3}&\t{if
$p=x^2+4y^2\e 1\mod 4$,}
\\\f {9}{98}R_1(p)-\f{9}{343}p
\mod {p^2}&\t{if $p\e 3\mod 4$ and $p\not=7$,}
\endcases
\\&\sum_{k=0}^{p-1}\f{k^3\b{2k}k^2\b{4k}{2k}}{648^k}
\\&\q\e \cases
\f{1436}{16807}x^2+\f{792}{16807}p-\f{1199p^2}{67228x^2}\mod
{p^3}&\t{if $p=x^2+4y^2\e 1\mod 4$,}
\\\f {216}{2401}R_1(p)-\f{1510}{16807}p
\mod {p^2}&\t{if $p\e 3\mod 4$ and $p\not=7$}\endcases
\endalign$$
and
$$\sum_{k=0}^{p-1}(7k+1)\f{\b{2k}k^2\b{4k}{2k}}{648^k}
\e (-1)^{\f{p-1}2}p-\f{745}{447}p^3E_{p-3}\mod {p^4}\q\t{for
$p\not=149$}.$$
\endpro

\pro{Conjecture 2.7} Let $p>3$ be a prime. Then
$$\align &\Ls
p3\sum_{k=0}^{p-1}\f{k\b{2k}k\b{3k}{k}\b{6k}{3k}}{12^{3k}} \e \cases
-\f{5}{9}x^2+\f{5p}{18}+\f{5p^2}{72x^2}\mod {p^3}&\t{if
$p=x^2+4y^2\e 1\mod 4$,}
\\\f 1{12}R_1(p)\mod {p^2}
&\t{if $p\e 3\mod 4$,}
\endcases
\\&\Ls
p3\sum_{k=0}^{p-1}\f{k^2\b{2k}k\b{3k}{k}\b{6k}{3k}}{12^{3k}}\\&\q \e
\cases \f{25}{486}x^2-\f{25}{972}p-\f{35p^2}{1944x^2}\mod
{p^3}&\t{if $p=x^2+4y^2\e 1\mod 4$,}
\\-\f {23}{648}R_1(p)\mod {p^2}
&\t{if $p\e 3\mod 4$,}
\endcases
\\&\Ls
p3\sum_{k=0}^{p-1}\f{k^3\b{2k}k\b{3k}{k}\b{6k}{3k}}{12^{3k}}
\\&\q\e\cases \f{5}{2187}x^2-\f{5}{4374}p+\f{187p^2}{69984x^2}\mod
{p^3}&\t{if $p=x^2+4y^2\e 1\mod 4$,}
\\\f {197}{29160}R_1(p)\mod {p^2}
&\t{if $p\e 3\mod 4$.}
\endcases
\endalign$$
\endpro

\pro{Conjecture 2.8} Let $p$ be a prime with $p\not=2,3,11$. Then
$$\align
&\Ls {33}p\sum_{k=0}^{p-1}\f{k^2\b{2k}k\b{3k}{k}\b{6k}{3k}}{66^{3k}}
\\&\q\e \cases
\f 1{27783}\big(370x^2-\f{1060}{3}p-\f{25p^2}{6x^2}\big)\mod
{p^3}&\t{if $p=x^2+4y^2\e 1\mod 4$,}
\\ \f {121}{7938}R_1(p)
+\f{505}{83349}p\mod {p^2}&\t{if $p\e 3\mod 4$ and $p>7$},
\endcases
\\&\Ls {33}p\sum_{k=0}^{p-1}\f{k^3\b{2k}k\b{3k}{k}\b{6k}{3k}}{66^{3k}}
\\&\q\e \cases
\f
1{4084101}\big(\f{21100}9x^2-\f{1300}{9}p-\f{485p^2}{4x^2}\big)\mod
{p^3}&\t{if $p=x^2+4y^2\e 1\mod 4$,}
\\\f 1{1750329}\big(242R_1(p)
-\f{9250}{21}p\big)\mod {p^2}&\t{if $p\e 3\mod 4$ and $p>7$.}
\endcases
\endalign$$
\endpro
\par{\bf Remark 2.3} Let $p$ be a prime with $p\not=2,3,11$. In [S5], the author
conjectured that
$$\Ls {33}p\sum_{k=0}^{p-1}(63k+5)
\f{\b{2k}k\b{3k}{k}\b{6k}{3k}}{66^{3k}}\e 5p\mod {p^3}.$$ For the
conjectures on $\sum_{k=0}^{p-1}\f{\b{2k}k^2\b{4k}{2k}}{648^k},\
\sum_{k=0}^{p-1}\f{\b{2k}k\b{3k}{k}\b{6k}{3k}}{12^{3k}}$ and
$\sum_{k=0}^{p-1} \f{\b{2k}k\b{3k}{k}\b{6k}{3k}}{66^{3k}}$ modulo
$p^3$ see [S7, Conjecture 4.20]. \vskip0.2cm
 \pro{Conjecture 2.9}
Let $p>3$ be a prime. Then
$$\align &\sum_{k=0}^{p-1}\f{k^2\b{2k}k^3}{16^k}\e
\cases \f 49x^2-\f 59 p+\f{p^2}{36x^2}\mod {p^3} &\t{if
$p=x^2+3y^2\e 1\mod 3$,}
\\-\f 29R_3(p)-\f 13p\mod{p^2}
&\t{if $p\e 2\mod 3$,}
\endcases
\\&\sum_{k=0}^{p-1}\f{k^3\b{2k}k^3}{16^k}\e
\cases -\f 29x^2+\f 49 p-\f{7p^2}{72x^2}\mod {p^3} &\t{if
$p=x^2+3y^2\e 1\mod 3$,}
\\\f 49R_3(p)+\f 13p\mod{p^2}
&\t{if $p\e 2\mod 3$.}
\endcases
\endalign$$
\endpro
\par{\bf Remark 2.4} Let $p>3$ be a prime. In [S7], the author
conjectured that $$\sum_{k=0}^{p-1}\f{\b{2k}k^3}{16^k}\e\cases
4x^2-2p-\f{p^2}{4x^2}\mod {p^3}&\t{if $p=x^2+3y^2\e 1\mod 3$,}
\\-p^2\b{(p-1)/2}{(p-5)/6}^{-2}\mod {p^3}&\t{if $p\e 2\mod 3$.}
\endcases$$ In [S11], the author conjectured that
$$\sum_{k=0}^{(p-1)/2}\f{\b{2k}k^3}{16^k(k+1)} \e\cases
-16y^2+2p\mod {p^3}&\t{if $p=x^2+3y^2\e 1\mod 3$,}
\\-\f 43R_3(p)-\f 23p\mod {p^2}&\t{if $p\e 2\mod 3$.}
\endcases$$

 \pro{Conjecture 2.10} Let $p>3$ be a prime. Then
$$\align
&(-1)^{\f{p-1}2}\sum_{k=0}^{p-1}\f{k^2\b{2k}k^3}{256^k} \e\cases \f
19x^2-\f p{18}-\f{p^2}{72x^2} \mod {p^3}&\t{if $p=x^2+3y^2\e 1\mod
3$,}
\\-\f 29R_3(p)\mod {p^2}&\t{if $p\e 2\mod 3$,}
\endcases
\\&(-1)^{\f{p-1}2}\sum_{k=0}^{p-1}\f{k^3\b{2k}k^3}{256^k}
\e\cases \f 1{18}x^2+\f p{72}-\f{p^2}{144x^2} \mod {p^3}&\t{if
$p=x^2+3y^2\e 1\mod 3$,}
\\-\f 19R_3(p)+\f 1{24}p\mod {p^2}&\t{if $p\e 2\mod 3$.}
\endcases
\endalign$$
 \endpro

\pro{Conjecture 2.11} Let $p>3$ be a prime. Then
$$\align &
\\&\sum_{k=0}^{p-1}\f{k\b{2k}k^2\b{3k}{k}}{108^k}
\e\cases -\f 89x^2+\f 49p+\f{p^2}{9x^2}\mod {p^3}&\t{if
$p=x^2+3y^2\e 1\mod 3$,}
\\-\f 49R_3(p)\mod {p^2}&\t{if $p\e 2\mod 3$.}
\endcases
\\&\sum_{k=0}^{p-1}\f{k^2\b{2k}k^2\b{3k}{k}}{108^k}
\e\cases \f {32}{243}x^2-\f{16}{243}p-\f{17p^2}{486x^2}\mod
{p^3}&\t{if $p=x^2+3y^2\e 1\mod 3$,}
\\\f {52}{243}R_3(p)\mod {p^2}&\t{if $p\e 2\mod 3$.}
\endcases
\\&\sum_{k=0}^{p-1}\f{k^3\b{2k}k^2\b{3k}{k}}{108^k}
\e\cases \f 1{10935}(16x^2-8p+\f{113p^2}{2x^2})\mod {p^3}&\t{if
$p=x^2+3y^2\e 1\mod 3$,}
\\-\f {92}{2187}R_3(p)\mod {p^2}&\t{if $p\e 2\mod 3$.}
\endcases
\endalign$$
\endpro

\pro{Conjecture 2.12} Let $p>3$ be a prime. Then
$$\align&\sum_{k=0}^{p-1}\f{k^2\b{2k}k^2\b{4k}{2k}}{(-144)^k} \e\cases
\f {12}{125}x^2-\f {19}{125}p+\f{7p^2}{500x^2}\mod {p^3}&\t{if
$p=x^2+3y^2\e 1\mod 3$,}
\\\f 2{25}R_3(p)+\f {13}{125}p\mod {p^2}&\t{if $p\e 2\mod 3$,}
\endcases
\\&\sum_{k=0}^{p-1}\f{k^3\b{2k}k^2\b{4k}{2k}}{(-144)^k} \e\cases
\f {62}{3125}x^2+\f {6}{3125}p-\f{511p^2}{25000x^2}\mod {p^3}&\t{if
$p=x^2+3y^2\e 1\mod 3$,}
\\-\f {48}{625}R_3(p)-\f {37}{3125}p\mod {p^2}&\t{if $p\e 2\mod 3$}
\endcases
\endalign$$
and
$$\sum_{k=0}^{p-1}(5k+1)\f{\b{2k}k^2\b{4k}{2k}}{(-144)^k}
\e (-1)^{[\f p3]}p+\f 52p^3U_{p-3}\mod {p^4}.$$

\endpro

\pro{Conjecture 2.13} Let $p>5$ be a prime. Then
$$\align &\Ls p5\sum_{k=0}^{p-1}\f{k^2\b{2k}k\b{3k}{k}\b{6k}{3k}}{54000^k}
\\&\q\e\cases \f {236}{11979}x^2-\f{199p}{11979}p-\f{37p^2}{47916x^2}\mod
{p^3}&\t{if $p=x^2+3y^2\e 1\mod 3$,}
\\\f {50}{1089}R_3(p)+\f{9}{1331}p\mod {p^2}&\t{if $p\e 2\mod 3$ and $p\not=11$,}
\endcases
\\&\Ls p5\sum_{k=0}^{p-1}\f{k^3\b{2k}k\b{3k}{k}\b{6k}{3k}}{54000^k}
\\&\q\e\cases \f 1{1449459}(3594x^2+500p-\f{1533p^2}{8x^2})\mod
{p^3}&\t{if $p=x^2+3y^2\e 1\mod 3$,}
\\\f {100}{43923}R_3(p)-\f{2297}{1449459}p\mod {p^2}&\t{if $p\e 2\mod 3$ and
$p\not=11$.}
\endcases
\endalign$$
\endpro
\par{\bf Remark 2.5} Let $p>5$ be a prime. In [S5], the author
conjectured that
$$\sum_{k=0}^{p-1}(11k+1)\f{\b{2k}k\b{3k}{k}\b{6k}{3k}}{54000^k}
\e \Ls{-15}pp\mod {p^3}.$$ In [S7, Conjecture 4.19], the author
conjectured that for $p\e 1\mod 3$ and so $p=x^2+3y^2$,
$$\sum_{k=0}^{p-1}\f{\b{2k}k^2\b{4k}{2k}}
{(-144)^k}\e \Ls
p5\sum_{k=0}^{p-1}\f{\b{2k}k\b{3k}{k}\b{6k}{3k}}{54000^k} \e
4x^2-2p-\f{p^2}{4x^2}\mod {p^3},$$ and for $p\e 2\mod 3$,
$$\sum_{k=0}^{p-1}\f{\b{2k}k^2\b{4k}{2k}}
{(-144)^k}\e 25\Ls
p5\sum_{k=0}^{p-1}\f{\b{2k}k\b{3k}{k}\b{6k}{3k}}{54000^k} \e
p^2\b{(p-1)/2}{(p-5)/6}^{-2}\mod {p^3}.$$

 \pro{Conjecture
2.14} Let $p>3$ be a prime. Then
$$\align &\sum_{k=0}^{p-1}\f{k^2\b{2k}k^2\b{3k}{k}}{1458^k}
\e\cases \f
{56}{1125}x^2-\f{14}{1125}(2+(-1)^{\f{p-1}2})p-\f{7p^2}{2250x^2}
\mod {p^3}\\\qq\qq\qq\ \qq\qq\t{if $p=x^2+3y^2\e 1\mod 3$,}
\\\f {8}{75}R_3(p)-\f{14}{1125}(-1)^{\f{p-1}2}p\mod {p^2}
\q\t{if $3\mid p-2$ and $p\not=5$,}
\endcases
\\&\sum_{k=0}^{p-1}\f{k^3\b{2k}k^2\b{3k}{k}}{1458^k}
\e\cases \f
1{84375}(904x^2-(452-524(-1)^{\f{p-1}2})p-\f{113p^2}{2x^2} \mod
{p^3}\\\qq\qq\qq\ \qq\qq\t{if $p=x^2+3y^2\e 1\mod 3$,}
\\\f {8}{625}R_3(p)+\f{524}{84375}(-1)^{\f{p-1}2}p\mod {p^2}
\q\t{if $3\mid p-2$ and $p\not=5$}
\endcases
\endalign$$
and
$$\sum_{k=0}^{p-1}(15k+2)\f{\b{2k}k^2\b{3k}{k}}{1458^k}
\e 2(-1)^{\f{p-1}2}p-\f{10}9p^3E_{p-3}\mod {p^4}.$$
\endpro
\par{\bf Remark 2.6} Let $p>3$ be a prime. In [S5] and [S7], the
author conjectured that
$$\sum_{k=0}^{p-1}\f{\b{2k}k^2\b{3k}{k}}{1458^k}
\e\cases 4x^2-2p-\f{p^2}{4x^2}\mod {p^3}&\t{if $p=x^2+3y^2\e 1\mod
3$,}\\0\mod {p^3}&\t{if $p\e 2\mod 3$.}\endcases$$ In [S11], the
author conjectured that
$$\sum_{k=0}^{p-1}\f{\b{2k}k^2\b{3k}{k}}{1458^k(k+1)} \e\cases
2(-1)^{\f{p-1}2}p-p^2\pmod {p^3}&\t{if $p\e 1\mod 3$,}
\\-12R_3(p)+2(-1)^{\f{p-1}2}p\pmod {p^2}&\t{if $p\e 2\mod 3$.}
\endcases$$

\pro{Conjecture 2.15} Let $p$ be an odd prime. Then
$$\align
&(-1)^{\f{p-1}2}\sum_{k=0}^{p-1}\f{k^2\b{2k}k^3}{(-64)^k} \e\cases
\f 18x^2-\f 3{16}p+\f{p^2}{64x^2}\mod{p^3}\q\t{if $p=x^2+2y^2\e
1,3\mod 8$,}
\\ -\f 1{32}R_2(p)-\f p8\mod{p^2}\q\t{if $p\e 5,7\mod 8$,}
\endcases
\\&(-1)^{\f{p-1}2}\sum_{k=0}^{p-1}\f{k^3\b{2k}k^3}{(-64)^k}
\e\cases \f 1{32}x^2-\f p{64}-\f{5p^2}{256x^2}\mod{p^3}\q\t{if
$p=x^2+2y^2\e 1,3\mod 8$,}
\\ \f 3{128}R_2(p)\mod{p^2}\q\t{if $p\e 5,7\mod 8$.}
\endcases
\endalign$$
\endpro

\pro{Conjecture 2.16} Let $p$ be an odd prime. Then
$$\align &\sum_{k=0}^{p-1}\f{k\b{2k}k^2\b{4k}{2k}}{256^k} \e\cases
-\f 34x^2+\f 38p+\f{3p^2}{32x^2}\mod {p^3}\\\qq\qq\qq\qq\q\ \t{if
$p=x^2+2y^2\e 1,3\mod 8$,}
\\-\f 1{16}R_2(p) \mod {p^2}\q\t{if $p\e 5,7\mod 8$,}
\endcases
\\&\sum_{k=0}^{p-1}\f{k^2\b{2k}k^2\b{4k}{2k}}{256^k} \e\cases
\f 3{32}x^2-\f 3{64}p-\f{7p^2}{256x^2}\mod {p^3}\\\qq\qq\qq\qq\q\
\t{if $p=x^2+2y^2\e 1,3\mod 8$,}
\\\f {11}{384}R_2(p) \mod {p^2}\q\ \t{if $p\e 5,7\mod 8$,}
\endcases
\\&\sum_{k=0}^{p-1}\f{k^3\b{2k}k^2\b{4k}{2k}}{256^k} \e\cases
\f 3{1280}x^2-\f 3{2560}p+\f{41p^2}{10240x^2}\mod
{p^3}\\\qq\qq\qq\qq\qq\t{if $p=x^2+2y^2\e 1,3\mod 8$,}
\\-\f {17}{3072}R_2(p) \mod {p^2}\ \t{if $p\e 5,7\mod 8$.}
\endcases
\endalign$$
\endpro

\pro{Conjecture 2.17} Let $p$ be a prime with $p\not=2,7$. Then
 $$\align
& \sum_{k=0}^{p-1}\f{k^2\b{2k}k^2\b{4k}{2k}}{28^{4k}} \e\cases \f
1{32000}(363x^2-2p-359\f{1+\sls p3}2p-\f{p^2}{2x^2}) \mod
{p^3}\\\qq\qq\qq\qq\qq\q\t{if $p=x^2+2y^2\e 1,3\mod 8$,}
\\\f {147}{25600}R_2(p)-\f{359}{64000}\ls p3p
\mod {p^2}\q\t{if $p\e 5,7\mod 8$,}
\endcases
\\& \sum_{k=0}^{p-1}\f{k^3\b{2k}k^2\b{4k}{2k}}{28^{4k}}
\e\cases \f 1{51200000}(3963x^2-2902p+1841\f{1+\sls
p3}2p-\f{451p^2}{2x^2}) \mod {p^3}\\\qq\qq\qq\qq\qq\qq\q\ \t{if
$p=x^2+2y^2\e 1,3\mod 8$,}
\\\f {1}{40960000}(147R_2(p)+\f{3682}{5}\ls p3p)
\mod {p^2}\q\t{if $p\e 5,7\mod 8$.}
\endcases
\endalign$$
\endpro
\par{\bf Remark 2.7} Let $p$ be a prime with $p\not=2,3,7$. In [Su1],
Z.W. Sun made a conjecture equivalent to
$$\sum_{k=0}^{p-1}(40k+3)\f{\b{2k}k^2\b{4k}{2k}}{28^{4k}}
\e 3\Ls p3p-\f{15}{196}p^3U_{p-3}\mod {p^4}.$$ For the conjectures
on $\sum_{k=0}^{p-1}\f{\b{2k}k^2\b{4k}{2k}}{256^k}$ and
$\sum_{k=0}^{p-1}\f{\b{2k}k^2\b{4k}{2k}}{28^{4k}}$ modulo $p^3$, see
[S7, Conjecture 4.18].

 \pro{Conjecture 2.18} Let $p$ be an odd
prime. Then
$$\align
&\sum_{k=0}^{p-1}\f{k^2\b{2k}k^2\b{3k}{k}}{8^k} \\&\e \cases
\f{87}{250}x^2-\f {213}{500}p+\f{39p^2}{2000x^2}\mod {p^3}&\t{if
$p=x^2+2y^2\e 1,3\mod 8$,}
\\-\f 3{200}R_2(p)-\f{63}{250}p\mod {p^2}&
\t{if $p\e 5,7\mod 8$ and $p\not=5$,}
\endcases
\\&\sum_{k=0}^{p-1}\f{k^3\b{2k}k^2\b{3k}{k}}{8^k}\\& \e \cases
-\f{1347}{12500}x^2+\f {5703}{25000}p-\f{6009p^2}{100000x^2}\mod
{p^3}&\t{if $p=x^2+2y^2\e 1,3\mod 8$,}
\\\f {243}{10000}R_2(p)+\f{1089}{6250}p\mod {p^2}&
\t{if $p\e 5,7\mod 8$ and $p\not=5$.}
\endcases
\endalign$$\endpro

\par{\bf Remark 2.8} Let $p$ be an odd prime. In [Su1],
Z.W. Sun conjectured that
$$\sum_{k=0}^{p-1}(10k+3)\f{\b{2k}k^2\b{3k}{k}}{8^k}
\e 3p+\f{49}8p^3B_{p-3}\mod {p^4}.$$ In [S7], the author conjectured
that
$$\sum_{k=0}^{p-1}\f{\b{2k}k^2\b{3k}{k}}{8^k}
\e\cases 4x^2-2p-\f{p^2}{4x^2}\mod {p^3}&\t{if $p=x^2+2y^2\e 1,3\mod
8$,}
\\\f {11}9p^2\b{[p/4]}{[p/8]}^{-2}\mod {p^3}
&\t{if $p\e 5\mod 8$,}
\\-\f {11}2p^2\b{[p/4]}{[p/8]}^{-2}\mod {p^3}
&\t{if $p\e 7\mod 8$.}\endcases$$ In [S11], the author conjectured
that $$\sum_{k=0}^{p-1}\f{\b{2k}k^2\b{3k}k}{8^k(k+1)} \e\cases
-11y^2+2p-p^2\mod {p^3}& \t{if $p=x^2+2y^2\e 1,3\mod 8$,}
\\-\f{1}8R_2(p)-\f 34p\mod {p^2}&\t{if $p\e 5,7\mod
8$.}\endcases$$

\pro{Conjecture 2.19} Let $p$ be a prime with $p\not=2,5$. Then
$$\align
&\Ls{-5}p\sum_{k=0}^{p-1}\f{k^2\b{2k}k\b{3k}k\b{6k}{3k}}{20^{3k}}
\\&\q\e\cases \f {111}{2744}x^2-\f{93}{5488}p-\f{129}{21952x^2}p^2\mod
{p^3}&\t{if $p=x^2+2y^2\e 1,3\mod 8$,}
\\-\f {25}{1568}R_2(p) +\f{9}{2744}p \mod {p^2}
&\t{if $p\e 5,7\mod 8$,}
\endcases
\\&\Ls{-5}p\sum_{k=0}^{p-1}\f{k^3\b{2k}k\b{3k}k\b{6k}{3k}}{20^{3k}}
\\&\q\e\cases \f 1{537824}(9579x^2+\f{12165}{2}p-\f{10743p^2}{8x^2})\mod
{p^3}&\t{if $p=x^2+2y^2\e 1,3\mod 8$,}
\\-\f {2025}{307328}R_2(p) +\f{1359}{67228}p \mod {p^2}
&\t{if $p\e 5,7\mod 8$.}
\endcases
\endalign$$
\endpro
\par{\bf Remark 2.9} Let $p>5$ be a prime. In [S5] the author
conjectured that
$$\sum_{k=0}^{p-1}(28k+3)\f{\b{2k}k\b{3k}k\b{6k}{3k}}{20^{3k}}
\e 3\Ls{-5}pp\mod {p^3}.$$ In [S7] the author conjectured the
congruence for
$\sum_{k=0}^{p-1}\f{\b{2k}k\b{3k}k\b{6k}{3k}}{20^{3k}}$ modulo
$p^3$. \vskip0.2cm
\par For any odd prime $p$ let
$$R_7(p)=\sum_{k=0}^{(p-1)/2}\f{\b{2k}k^3}{k+1}.$$
In [S11] the author conjectured that for any prime $p\not=2,3,7$,
$$R_7(p)
\e \cases -44y^2+2p\mod {p^3}&\t{if $p=x^2+7y^2\e 1,2,4\mod 7$,}
\\-\f{1}7\b{[3p/7]}{[p/7]}^2\mod p&\t{if $p\e 3\mod 7$,}
\\-\f{16}7\b{[3p/7]}{[p/7]}^2\mod p&\t{if $p\e 5\mod 7$,}
\\-\f{4}7\b{[3p/7]}{[p/7]}^2\mod p&\t{if $p\e 6\mod 7$.}
\endcases$$

\pro{Conjecture 2.20} Let $p$ be a prime with $p\not=2,3,7$. Then
$$\align&\sum_{k=0}^{p-1}k^2\b{2k}k^3\\&\e
\cases \f{736x^2}{1323}-\f{272p}{441}+\f{20p^2}{1323x^2}\mod {p^3}
&\t{if $p=x^2+7y^2\e 1,2,4\mod 7$,}
\\\f 8{63}R_7(p)-\f{256}{1323}p\mod {p^2}
&\t{if $p\e 3,5,6\mod 7$,}
\endcases
\\&\sum_{k=0}^{p-1}k^3\b{2k}k^3\\&\e \cases
-\f{5408x^2}{27783}+\f{2992p}{9261}-\f{1774p^2}{27783x^2}\mod {p^3}
&\t{if $p=x^2+7y^2\e 1,2,4\mod 7$,}
\\-\f {256}{1323}R_7(p)+\f{128}{27783}p\mod {p^2}
&\t{if $p\e 3,5,6\mod 7$.}
\endcases\endalign$$
\endpro

\pro{Conjecture 2.21} Let $p$ be a prime with $p\not=2,3,7$. Then
$$\align
&(-1)^{\f{p-1}2}\sum_{k=0}^{(p-1)/2}k^2\f{\b{2k}k^3}{4096^k}
\\&\q\e \cases \f{43x^2}{1323}-\f{13}{441}p-\f{p^2}{1323x^2}\mod {p^3}&\t{if $p=x^2+7y^2\e 1,2,4\mod 7$,}
\\\f 8{63}R_7(p)+\f{349}{2646}p\mod {p^2}&\t{if $p\e 3,5,6\mod 7$,}
\endcases
\\&(-1)^{\f{p-1}2}\sum_{k=0}^{(p-1)/2}k^3\f{\b{2k}k^3}{4096^k}
\\&\q\e \cases \f{169x^2}{55566}-\f{31}{74088}p-\f{71p^2}{444528x^2}\mod {p^3}&\t{if $p=x^2+7y^2\e 1,2,4\mod 7$,}
\\\f 4{1323}R_7(p)+\f{1013}{222264}p\mod {p^2}&\t{if $p\e 3,5,6\mod 7$.}
\endcases\endalign$$
\endpro

\pro{Conjecture 2.22} Let $p$ be a prime with $p>7$. Then
$$\align
&\sum_{k=0}^{p-1}k^2\f{\b{2k}k^2\b{4k}{2k}}{81^k} \e \cases
\f{1376}{6125}x^2-\f{2032}{6125}p+\f{164p^2}{6125x^2}\mod
{p^3}\\\qq\qq\qq\qq\t{if $p=x^2+7y^2\e 1,2,4\mod 7$,}
\\\f{36}{175}R_7(p)+\f {96}{6125}p\mod {p^2}\q\t{if $p\e 3,5,6\mod 7$,}
\endcases
\\&\sum_{k=0}^{p-1}k^3\f{\b{2k}k^2\b{4k}{2k}}{81^k}
\e \cases \f
1{1071785}(-130784x^2+335088p-87826\f{p^2}{x^2})\mod{p^3}\\\qq\qq\qq\qq\qq\qq\t{if
$p=x^2+7y^2\e 1,2,4\mod 7$,}
\\-\f{13824}{30625}R_7(p)-\f{283264}{1071875}p\mod {p^2}\q\t{if $p\e 3,5,6\mod 7$.}
\endcases\endalign$$
\endpro
\par{\bf Remark 2.10} Let $p>3$ be a prime. In [Su1], Z.W. Sun
conjectured the congruence for
$\sum_{k=0}^{p-1}\f{\b{2k}k^2\b{4k}{2k}}{81^k}$ modulo $p^2$, and
$$\sum_{k=0}^{p-1}(35k+8)\f{\b{2k}k^2\b{4k}{2k}}{81^k}\e
8p+\f{416}{27}p^3B_{p-3}\mod {p^4}.$$ In [S7], the author
conjectured the congruence for
$\sum_{k=0}^{p-1}\f{\b{2k}k^2\b{4k}{2k}}{81^k}$ modulo $p^3$.

\pro{Conjecture 2.23} Let $p$ be a prime with $p\not=2,3,5,7,13$.
Then
$$\align
& \sum_{k=0}^{p-1}k^2\f{\b{2k}k^2\b{4k}{2k}} {(-3969)^k} \e \cases
\f 1{274625}(6176x^2-8272p+\f{524p^2}{x^2}) \mod
{p^3}\\\qq\qq\qq\qq\qq\t{if $p=x^2+7y^2\e 1,2,4\mod 7$,}
\\-\f{504}{4225}R_7(p)-\f{32256}{274625}p\mod {p^2}\q\t{if $p\e 3,5,6\mod 7$,}
\endcases
\\& \sum_{k=0}^{p-1}k^3\f{\b{2k}k^2\b{4k}{2k}} {(-3969)^k} \e \cases
\f 1{65^5}(-4940384x^2-1487952p-\f{40666p^2}{x^2})\mod
{p^3}\\\qq\qq\qq\qq\qq\t{if $p=x^2+7y^2\e 1,2,4\mod 7$,}
\\\f{193536}{65^4}R_7(p)+\f{18335104}{65^5}p\mod {p^2}\q\t{if $p\e 3,5,6\mod 7$}
\endcases
 \endalign$$
 and
 $$\sum_{k=0}^{p-1}(65k+8)\f{\b{2k}k^2\b{4k}{2k}} {(-3969)^k}
 \e 8\Ls p7p\mod{p^3}.$$
\endpro

\pro{Conjecture 2.24} Let $p$ be a prime with $p>7$. Then
$$\align&\Ls {-15}p\sum_{k=0}^{p-1}k^2\f{\b{2k}k\b{3k}k\b{6k}{3k}}
{(-15)^{3k}} \e \cases \f
1{1323}(32x^2-\f{592}9p+\f{76p^2}{9x^2})\mod
{p^3}\\\qq\qq\qq\qq\t{if $p=x^2+7y^2\e 1,2,4\mod 7$,}
\\\f {50}{567}R_7(p)+\f{752}{11907}p\mod {p^2}\q\t{if $p\e 3,5,6\mod 7$,}
\endcases
\\&\Ls {-15}p\sum_{k=0}^{p-1}k^3\f{\b{2k}k\b{3k}k\b{6k}{3k}}
{(-15)^{3k}} \e \cases \f
1{750141}(2656x^2-7600p-\f{3174p^2}{x^2})\mod
{p^3}\\\qq\qq\qq\qq\qq\t{if $p=x^2+7y^2\e 1,2,4\mod 7$,}
\\-\f {1600}{35721}R_7(p)-\f{44672}{750141}p\mod {p^2}\ \t{if $p\e 3,5,6\mod 7$}
\endcases
\endalign$$
and
$$\sum_{k=0}^{p-1}(63k+8)\f{\b{2k}k\b{3k}k\b{6k}{3k}}
{(-15)^{3k}} \e 8\Ls{-15}pp\mod {p^3}.$$

\endpro

\pro{Conjecture 2.25} Let $p$ be a prime with $p\not=2,3,11$. Then
$$\align&\sum_{k=0}^{p-1}k^2\f{\b{2k}k^2\b{3k}{k}}{64^k}
\\&\e\cases \f{12}{121}x^2-\f{147}{242}p+\f{51p^2}{242x^2}\mod {p^3}
&\t{if $\ls p{11}=1$ and so $4p=x^2+11y^2$,}
\\\f 3{11}\sum_{k=0}^{p-1}\f{\b{2k}k^2\b{3k}{k}}{64^k(k+1)}
-\f{81}{242}p\mod {p^2}&\t{if $\ls p{11}=-1$,}\endcases
\\&\sum_{k=0}^{p-1}k^3\f{\b{2k}k^2\b{3k}{k}}{64^k}
\\&\e\cases \f 1{2662}(-411x^2+\f{7089}{2}p-\f{5253p^2}{2x^2})\mod {p^3}
&\t{if $\ls p{11}=1$ and so $4p=x^2+11y^2$,}
\\-\f {243}{242}\sum_{k=0}^{p-1}\f{\b{2k}k^2\b{3k}{k}}{64^k(k+1)}
+\f{3987}{5324}p\mod {p^2}&\t{if $\ls p{11}=-1$.}\endcases
\endalign$$
\endpro

\pro{Conjecture 2.26} Let $p$ be a prime with $p\not=2,3,11$. Then
$$\align&\Ls{-2}p\sum_{k=0}^{p-1}k^2\f{\b{2k}k\b{3k}{k}\b{6k}{3k}}{(-32)^{3k}}
\e\cases
\f{615}{166012}x^2-\f{1605}{83006}p+\f{375p^2}{83006x^2}\mod {p^3}
\\\qq\qq\q\qq\t{if $\ls p{11}=1$ and so $4p=x^2+11y^2$,}
\\\f{16}{539}\sum_{k=0}^{p-1}\f{\b{2k}k^2\b{3k}{k}}{64^k(k+1)}-\f{159}{41503}p
\mod {p^2}\\\qq\qq\qq\q\t{if $\ls p{11}=-1$,}
\endcases
\\&\Ls{-2}p\sum_{k=0}^{p-1}k^3\f{\b{2k}k\b{3k}{k}\b{6k}{3k}}{(-32)^{3k}}
\e\cases \f 1{7^5\cdot
22^3}\big(-114495x^2-168495p+\f{3930p^2}{x^2}\big)\mod {p^3}
\\\qq\qq\q\qq\t{if $\ls p{11}=1$ and so $4p=x^2+11y^2$,}
\\\f 1{7^5\cdot 22^3}\big(-399168\sum_{k=0}^{p-1}\f{\b{2k}k^2\b{3k}{k}}{64^k(k+1)}
-506349p\big) \mod {p^2}\\\qq\qq\q\qq\t{if $\ls p{11}=-1$.}
\endcases
\endalign$$
\endpro

\par{\bf Remark 2.11} Let $p$ be an odd prime with $\ls p{11}=-1$. In [S11],
the author conjectured that
$$ \sum_{k=0}^{p-1}\f{\b{2k}k^2\b{3k}{k}}
{64^k(k+1)}\e \cases -\f{50}{11}R_{11}(p)\mod p&\t{if $p\e 2\mod
{11}$,}
\\-\f{32}{11}R_{11}(p)\mod p&\t{if $p\e 6\mod{11}$,}
\\-\f{2}{11}R_{11}(p)\mod p&\t{if $p\e 7\mod{11}$,}
\\-\f{72}{11}R_{11}(p)\mod p&\t{if $p\e 8\mod{11}$,}
\\-\f{18}{11}R_{11}(p)\mod p&\t{if $p\e 10\mod{11}$,}
\endcases$$
where $$R_{11}(p)=\b{[\f{3p}{11}]}{[\f
p{11}]}^2\b{[\f{6p}{11}]}{[\f{3p}{11}]}^2\cdot\b{[\f{4p}{11}]}{[\f{2p}{11}]}^{-2}.$$
In [Su1], Z.W. Sun conjectured that for any prime $p>3$,
$$\sum_{k=0}^{p-1}(11k+3)\f{\b{2k}k^2\b{3k}{k}}
{64^k}\e 3p+\f 72p^3B_{p-3}\mod {p^4}.$$ Z.W. Sun also conjectured
that
$$\sum_{k=0}^{p-1}(154k+15)\f{\b{2k}k\b{3k}{k}\b{6k}{3k}}{(-32)^{3k}}
\e 15\Ls{-2}pp\mod {p^3}.$$
 For the congruences on
$\sum_{k=0}^{p-1}\f{\b{2k}k^2\b{3k}{k}} {64^k}$ and
$\sum_{k=0}^{p-1}\f{\b{2k}k\b{3k}{k}\b{6k}{3k}}{(-32)^{3k}}$ modulo
$p^3$, see [S7, Conjecture 4.22].

 \pro{Conjecture 2.27} Let $p$ be a prime with $p\not=2,3,19$.
Then
$$\align&\Ls{-6}p\sum_{k=0}^{p-1}k^2\f{\b{2k}k\b{3k}{k}\b{6k}{3k}}{(-96)^{3k}}
\\&\e\cases \f 1{342^2}\big(305x^2-\f{3730}{3}p+\f{70p^2}{3x^2}\big)\mod {p^3}
\\\qq\qq\qq\qq\qq\qq\qq\t{if $\ls p{19}=1$ and so $4p=x^2+19y^2$,}
\\\f 1{342^2}\big(\f{95}6\ls{-6}p\sum_{k=0}^{p-1}\f{\b{2k}k\b{3k}{k}\b{6k}{3k}}{(-96)^{3k}(k+1)}
-\f{25}{3}p\big)\mod {p^2}\q\t{if $\ls p{19}=-1$,}\endcases
\\&\Ls{-6}p\sum_{k=0}^{p-1}k^3\f{\b{2k}k\b{3k}{k}\b{6k}{3k}}{(-96)^{3k}}
\\&\e\cases \f 1{342^3}\big(-\f{6235}3x^2+\f{3445}{3}p+\f{1750p^2}{x^2}\big)\mod {p^3}
\\\qq\qq\qq\qq\qq\qq\qq\qq\t{if $\ls p{19}=1$ and so $4p=x^2+19y^2$,}
\\\f 1{342^3}\big(-\f{95}6\ls{-6}p\sum_{k=0}^{p-1}\f{\b{2k}k\b{3k}{k}\b{6k}{3k}}{(-96)^{3k}(k+1)}
-\f{10900}{3}p\big)\mod {p^2}\q\t{if $\ls p{19}=-1$.}\endcases
\endalign$$
\endpro
\par{\bf Remark 2.12} Let $p>3$ be a prime. For the congruences on
$\sum_{k=0}^{p-1}\f{\b{2k}k\b{3k}{k}\b{6k}{3k}}{(-96)^{3k}}$ modulo
$p^3$ and
$\sum_{k=0}^{p-1}\f{\b{2k}k\b{3k}{k}\b{6k}{3k}}{(-96)^{3k}(k+1)}$
modulo $p$ see [S7,S10]. In [Su1], Z.W. Sun conjectured that
$$\sum_{k=0}^{p-1}(342k+25)\f{\b{2k}k\b{3k}{k}\b{6k}{3k}}{(-96)^{3k}}
\e 25\Ls{-6}pp\mod {p^3}.$$

\pro{Conjecture 2.28} Let $p$ be a prime with $p\not=2,3,7$. Then
$$\align&\sum_{k=0}^{p-1}k^2\f{\b{2k}k^2\b{4k}{2k}}{(-12288)^k}
\e\cases \f{27}{1372}x^2-\f{123}{5488}p+\f{15p^2}{21952x^2}\mod
{p^3} \\\qq\qq\qq\t{if $12\mid p-1$ and so $p=x^2+9y^2$,}
\\-\f{27}{2744}x^2+\f{123}{5488}p-\f{15p^2}{10976x^2}\mod {p^3}
\\\qq\qq\qq\t{if $12\mid p-5$ and so $2p=x^2+9y^2$,}
\\\f 3{1568}\sum_{k=0}^{p-1}\f{\b{2k}k^2\b{4k}{2k}}{(-12288)^k(k+1)}
-\f 3{2744}\ls p3p\mod {p^2}\\\qq\qq\qq\t{if $p\e 3\mod {4}$,}
\endcases
\\&\sum_{k=0}^{p-1}k^3\f{\b{2k}k^2\b{4k}{2k}}{(-12288)^k}
\e\cases
-\f{603}{268912}x^2+\f{3}{1075648}p+\f{351p^2}{4302592x^2}\mod {p^3}
\\\qq\qq\qq\t{if $12\mid p-1$ and so $p=x^2+9y^2$,}
\\\f 1{537824}(603x^2-\f 32p-\f{351p^2}{4x^2})\mod {p^3}
\\\qq\qq\qq\t{if $12\mid p-5$ and so $2p=x^2+9y^2$,}
\\-\f 9{153664}\sum_{k=0}^{p-1}\f{\b{2k}k^2\b{4k}{2k}}{(-12288)^k(k+1)}
-\f {1581}{1075648}\ls p3p\mod {p^2}\\\qq\qq\qq\t{if $p\e 3\mod
{4}$.}
\endcases
\endalign$$
\endpro
\par {\bf Remark 2.13} Let $p>3$ be a prime. In [Su1], Z.W. Sun
conjectured the congruence for
$\sum_{k=0}^{p-1}\f{\b{2k}k^2\b{4k}{2k}}{(-12288)^k}$ modulo $p^2$
and made a conjecture equivalent to
$$\sum_{k=0}^{p-1}(28k+3)\f{\b{2k}k^2\b{4k}{2k}}{(-12288)^k}
\e 3\Ls p3p+\f 54p^3U_{p-3}\mod {p^4}.$$ In [S7], the author
conjectured that
$$\sum_{k=0}^{p-1}\f{\b{2k}k^2\b{4k}{2k}}{(-12288)^k}
\e\cases 4x^2-2p-\f{p^2}{4x^2}\mod {p^3}&\t{if $12\mid p-1$ and so
$p=x^2+9y^2$,}
\\-2x^2+2p+\f{p^2}{2x^2}\mod {p^3}&\t{if $12\mid p-5$ and so
$2p=x^2+9y^2$,}
\\-\f 16p^2\b{[p/3]}{[p/12]}^{-2}\mod {p^3}&\t{if $p\e 7\mod {12}$,}
\\\f 1{12}p^2\b{[p/3]}{[p/12]}^{-2}\mod {p^3}&\t{if $p\e 11\mod {12}$.}
\endcases$$

\pro{Conjecture 2.29} Let $p$ be a prime with $p\not=2,5$.
 Then
$$\align&\sum_{k=0}^{p-1}k^2\f{\b{2k}k^2\b{4k}{2k}}{(-1024)^k}
\e\cases \f 3{100}x^2-\f{21}{400}p+\f{9p^2}{1600x^2}\mod
{p^3}\\\qq\qq\qq\qq\q\t{if $p\e 1,9\mod{20}$ and so
 $p=x^2+5y^2$,}
\\-\f 3{200}x^2+\f{21}{400}p-\f{9p^2}{800x^2}\mod {p^3}\\\qq\qq\qq\qq\q\t{if
$p\e 3,7\mod{20}$ and so
 $2p=x^2+5y^2$,}
 \\\f 3{160}\sum_{k=0}^{p-1}\f{\b{2k}k^2\b{4k}{2k}}{(-1024)^k(k+1)}
 -\f 3{200}(-1)^{\f{p-1}2}p\mod {p^2}\\\qq\qq\qq\qq\q \t{if $p\e 11,13,17,19\mod
 {20}$,}
 \endcases
\\&\sum_{k=0}^{p-1}k^3\f{\b{2k}k^2\b{4k}{2k}}{(-1024)^k}
\e\cases -\f 3{2000}x^2-\f{69}{8000}p-\f{69p^2}{32000x^2}\mod
{p^3}\\\qq\qq\qq\qq\q\t{if $p\e 1,9\mod{20}$ and so
 $p=x^2+5y^2$,}
\\\f 3{4000}x^2+\f{69}{8000}p+\f{69p^2}{16000x^2}\mod {p^3}\\\qq\qq\qq\qq\q\t{if
$p\e 3,7\mod{20}$ and so
 $2p=x^2+5y^2$,}
 \\-\f 9{1600}\sum_{k=0}^{p-1}\f{\b{2k}k^2\b{4k}{2k}}{(-1024)^k(k+1)}
 -\f {129}{8000}(-1)^{\f{p-1}2}p\mod {p^2}\\\qq\qq\qq\qq\q \t{if $p\e 11,13,17,19\mod
 {20}$.}
 \endcases
\endalign$$
\endpro
\par{\bf Remark 2.14} For any prime $p\not=2,5$, the congruence for
$\sum_{k=0}^{p-1}\b{2k}k^2\b{4k}{2k}(-1024)^{-k}$ modulo $p^2$ was
first conjectured by Z.W. Sun in [Su1]. Let
$R_{20}(p)=\b{\f{p-1}2}{[p/20]}\b{\f{p-1}2}{[3p/20]}$. In [S11], the
author conjectured that
$$\sum_{k=0}^{p-1}\f{\b{2k}k^2\b{4k}{2k}}{(-1024)^k(k+1)}
\e\cases \f 85R_{20}(p)\mod p& \t{if $p\e 1,3,7,9\mod {20}$,}
\\\f{4}{5}R_{20}(p)\mod p& \t{if $p\e 11\mod {20}$,}
\\\f{36}{5}R_{20}(p)\mod p& \t{if $p\e 13\mod {20}$,}
\\\f{28}{15}R_{20}(p)\mod p& \t{if $p\e 17\mod {20}$,}
\\\f{84}{5}R_{20}(p)\mod p& \t{if $p\e 19\mod {20}$.}
\endcases$$

\pro{Conjecture 2.30} Let $p>3$ be a prime. Then
$$\align&\sum_{k=0}^{p-1}k^2\f{\b{2k}k^2\b{3k}k}{216^k}
\\&\e\cases \f 16x^2-\f 1{36}p-\f{5p^2}{144x^2}\mod {p^3} & \t{if
$p=x^2+6y^2\e 1,7\mod {24}$,}
\\\f 13x^2-\f 5{36}p-\f{5p^2}{288x^2}\mod {p^3}
&\t{if $p=2x^2+3y^2\e 5,11\mod {24}$,}
\\-\f 19\sum_{k=0}^{p-1}\f{\b{2k}k^2\b{3k}k}{216^k(k+1)}
+\f 1{12}\ls p3p\mod {p^2}&\t{if $p\e 13,17,19,23\mod{24}$,}
\endcases
\\&
\sum_{k=0}^{p-1}k^3\f{\b{2k}k^2\b{3k}k}{216^k}
\\&\e\cases \f {23}{108}x^2+\f 7{72}p-\f{43p^2}{864x^2}\mod {p^3} & \t{if
$p=x^2+6y^2\e 1,7\mod {24}$,}
\\\f {23}{54}x^2-\f {67}{216}p-\f{43p^2}{1728x^2}\mod {p^3}
&\t{if $p=2x^2+3y^2\e 5,11\mod {24}$,}
\\-\f 16\sum_{k=0}^{p-1}\f{\b{2k}k^2\b{3k}k}{216^k(k+1)}
+\f {53}{216}\ls p3p\mod {p^2}&\t{if $p\e 13,17,19,23\mod{24}$.}
\endcases
\endalign$$
\endpro
\par{\bf Remark 2.15} Let $p>3$ be a prime and
$R_{24}(p)=\b{\f{p-1}2}{[\f p{24}]}\b{\f{p-1}2}{[\f{5p}{24}]}$. In
[S11], the author conjectured that
$$\sum_{k=0}^{p-1}\f{\b{2k}k^2\b{3k}k}{216^k(k+1)}\e\cases
\f{7}{8}R_{24}(p) \mod p& \t{if $p\e 1,11\mod {24}$,}
\\-\f{7}{8}R_{24}(p) \mod p& \t{if
$p\e 5,7\mod {24}$,}
\\ -\f 5{8}R_{24}(p)
\mod p& \t{if $p\e 13\mod {24}$,}
\\\f 5{8}R_{24}(p) \mod
p& \t{if $p\e 17\mod {24}$,}
\\\f {77}{8}R_{24}(p) \mod
p& \t{if $p\e 19\mod {24}$,}
\\-\f {77}{8}R_{24}(p) \mod
p& \t{if $p\e 23\mod {24}$.}
\endcases$$

\pro{Conjecture 2.31} Let $p>3$ be a prime. Then
$$\align&\Ls p3\sum_{k=0}^{p-1}k^2\f{\b{2k}k^2\b{4k}{2k}}{48^{2k}}
\\&\e\cases \f 3{64}x^2-\f 1{32}p-\f{p^2}{256x^2}\mod {p^3} & \t{if
$p=x^2+6y^2\e 1,7\mod {24}$,}
\\\f 3{32}x^2-\f 1{32}p-\f{p^2}{512x^2}\mod {p^3}
&\t{if $p=2x^2+3y^2\e 5,11\mod {24}$,}
\\\f 1{32}\sum_{k=0}^{p-1}\f{\b{2k}k^2\b{3k}k}{216^k(k+1)}
-\f {1+\sls p3}{128}p\mod {p^2}&\t{if $p\e 13,17,19,23\mod{24}$,}
\endcases
\\&\Ls p3
\sum_{k=0}^{p-1}k^3\f{\b{2k}k^2\b{4k}{2k}}{48^{2k}}
\\&\e\cases \f {13}{1024}x^2+\f 3{1024}p-\f{p^2}{1024x^2}\mod {p^3} & \t{if
$p=x^2+6y^2\e 1,7\mod {24}$,}
\\\f {13}{512}x^2+\f {3}{1024}p-\f{p^2}{2048x^2}\mod {p^3}
&\t{if $p=2x^2+3y^2\e 5,11\mod {24}$,}
\\\f 3{512}\sum_{k=0}^{p-1}\f{\b{2k}k^2\b{3k}k}{216^k(k+1)}
+\f {1}{128}p\mod {p^2}&\t{if $p\e 13,19\mod{24}$,}
\\\f 3{512}\sum_{k=0}^{p-1}\f{\b{2k}k^2\b{3k}k}{216^k(k+1)}
+\f {11}{1024}p\mod {p^2}&\t{if $p\e 17,23\mod{24}$.}
\endcases
\endalign$$
\endpro

\par{\bf Remark 2.16} Let $p$ be a prime with $p>3$. In
[S7, Conjecture 4.24], the author conjectured the congruences for
$\sum_{k=0}^{p-1}\f{\b{2k}k^2\b{3k}k}{216^k}$ and
$\sum_{k=0}^{p-1}\f{\b{2k}k^2\b{4k}{2k}}{48^{2k}}$ modulo $p^3$ in
the case $\sls{-6}p=1$. The corresponding congruences modulo $p^2$
were conjectured by Z.W. Sun in [Su1]. In 2019, Guo and Zudilin
[GuoZ] proved Z.W. Sun's conjecture
$$\sum_{k=0}^{p-1}(8k+1)\f{\b{2k}k^2\b{4k}{2k}}{48^{2k}}
\e \Ls p3p\mod {p^3}.$$

\pro{Conjecture 2.32} Let $p>5$ be a prime. Then
$$\align &\sum_{k=0}^{p-1}k^2\f{\b{2k}k^2\b{3k}{k}} {(-27)^k}
\\&\e\cases \f{16}{75}x^2-\f{64}{225}p+\f{4p^2}{225x^2}\mod {p^3}
&\t{if $p=x^2+15y^2\e 1,19\mod{30}$,}
\\-\f{16}{25}x^2+\f{64}{225}p-\f{4p^2}{675x^2}\mod {p^3}
&\t{if $p=3x^2+5y^2\e 17,23\mod {30},$}
\\\f 4{45}\sum_{k=0}^{p-1}\f{\b{2k}k^2\b{3k}{k}} {(-27)^k(k+1)}
-\f 8{75}\ls p3p\mod {p^2} &\t{if $p\e 7,11,13,29\mod {30}$,}
\endcases
\\&\sum_{k=0}^{p-1}k^3\f{\b{2k}k^2\b{3k}{k}} {(-27)^k}
\\&\e\cases \f{16}{3375}x^2+\f{64}{1125}p-\f{127p^2}{3375x^2}\mod {p^3}
&\t{if $p=x^2+15y^2\e 1,19\mod{30}$,}
\\-\f{16}{1125}x^2-\f{64}{1125}p+\f{127p^2}{10125x^2}\mod {p^3}
&\t{if $p=3x^2+5y^2\e 17,23\mod {30},$}
\\-\f 8{75}\sum_{k=0}^{p-1}\f{\b{2k}k^2\b{3k}{k}} {(-27)^k(k+1)}
-\f {88}{3375}\ls p3p\mod {p^2} &\t{if $p\e 7,11,13,29\mod {30}$.}
\endcases
\endalign$$
\endpro

\pro{Conjecture 2.33} Let $p>5$ be a prime. Then
$$\align &\sum_{k=0}^{p-1}k^2\f{\b{2k}k^2\b{3k}{k}} {15^{3k}}
\\&\e\cases \f{48}{1331}x^2-\f{368}{11979}p-\f{16p^2}{11979x^2}\mod {p^3}
&\t{if $p=x^2+15y^2\e 1,19\mod{30}$,}
\\\f 1{1331}(-144x^2+\f{368}{9}p+\f{16p^2}{27x^2})\mod {p^3}
&\t{if $p=3x^2+5y^2\e 17,23\mod {30},$}
\\\f {80}{1089}\sum_{k=0}^{p-1}\f{\b{2k}k^2\b{3k}{k}} {(-27)^k(k+1)}
+\f {184}{3993}\ls p3p\mod {p^2} &\t{if $p\e 7,11,13,29\mod {30}$,}
\endcases
\\&\sum_{k=0}^{p-1}k^3\f{\b{2k}k^2\b{3k}{k}} {15^{3k}}
\\&\e\cases \f 1{27\cdot 11^5}(21328x^2-816p-\f{1135p^2}{x^2})\mod {p^3}
&\t{if $p=x^2+15y^2\e 1,19\mod{30}$,}
\\\f{1}{27\cdot 11^5}(-63984x^2+816p+\f{1135p^2}{3x^2})\mod {p^3}
&\t{if $p=3x^2+5y^2\e 17,23\mod {30},$}
\\\f {160}{3\cdot 11^4}\sum_{k=0}^{p-1}\f{\b{2k}k^2\b{3k}{k}} {(-27)^k(k+1)}
+\f {22520}{27\cdot 11^5}\ls p3p\mod {p^2} &\t{if $p\e
7,11,13,29\mod {30}$}
\endcases
\endalign$$
and
$$\sum_{k=0}^{p-1}(33k+4)\f{\b{2k}k^2\b{3k}{k}} {15^{3k}}
\e 4\Ls p3p-\f{52}{25}p^3U_{p-3}\mod {p^4}.$$
\endpro
\par{\bf Remark 2.17} In [S11], the author conjectured that for any
prime $p\e 7,11,13,29\mod {30}$,
$$\sum_{k=0}^{p-1}\f{\b{2k}k^2\b{3k}{k}}
{(-27)^k(k+1)}\e\cases \f 25\cdot 5^{-[\f
p3]}\b{[p/3]}{[p/15]}^2\mod p&\t{if $p\e 7\mod {30}$,}
\\\f  1{10}\cdot 5^{-[\f p3]}\b{[p/3]}{[p/15]}^2\mod p&\t{if $p\e 11\mod
{30}$,}
\\\f{32}{5}\cdot 5^{-[\f p3]}\b{[p/3]}{[p/15]}^2\mod p&\t{if $p\e 13\mod
{30}$,}
\\\f{8}{5}\cdot 5^{-[\f p3]}\b{[p/3]}{[p/15]}^2\mod p&\t{if $p\e 29\mod
{30}$.}
\endcases$$
The congruence for $\sum_{k=0}^{p-1}\f{\b{2k}k^2\b{3k}{k}}
{(-27)^k}\mod {p^2}$ was conjectured by Z.W. Sun in [Su1], and the
congruence for $\sum_{k=0}^{p-1}\f{\b{2k}k^2\b{3k}{k}} {(-27)^k}\mod
{p^3}$ was conjectured by the author in [S7]. In [Su1], Z.W. Sun
made a conjecture equivalent to
$$\sum_{k=0}^{p-1}(15k+4)\f{\b{2k}k^2\b{3k}{k}} {(-27)^k}\e 4\Ls
p3p+8p^3U_{p-3} \mod {p^4}.$$

\pro{Conjecture 2.34} Let $p$ be a prime with $p>5$. Then
$$\align&\sum_{k=0}^{p-1}k^2\f{\b{2k}k^2\b{4k}{2k}}{12^{4k}}
\e \cases\f{71}{3200}x^2-\f{131}{6400}p-\f{11p^2}{25600x^2}\mod
{p^3} \\\qq\qq\qq\qq\t{if $p=x^2+10y^2\e 1,9,11,19\mod {40}$,}
\\-\f{71}{1600}x^2+\f{131}{6400}p+\f{11p^2}{51200x^2}\mod {p^3}
\\\qq\qq\qq\qq\t{if $p=2x^2+5y^2\e 7,13,23,37\mod {40}$,}
\\-\f 3{2560}\sum_{k=0}^{p-1}\f{\b{2k}k^2\b{4k}{2k}}{12^{4k}(k+1)}
-\f 1{1600}\Ls p5p\mod{p^2}\ \t{if $\ls{-10}p=-1$,}
\endcases
\\&\sum_{k=0}^{p-1}k^3\f{\b{2k}k^2\b{4k}{2k}}{12^{4k}} \e
\cases\f 1{80^3}(853x^2-\f{153}{2}p-\f{353p^2}{8x^2})\mod {p^3}
\\\qq\qq\qq\qq\t{if $p=x^2+10y^2\e 1,9,11,19\mod {40}$,}
\\\f{1}{256000}(-853x^2+\f{153}{4}p+\f{353p^2}{32x^2})\mod {p^3}
\\\qq\qq\qq\qq\t{if $p=2x^2+5y^2\e 7,13,23,37\mod {40}$,}
\\ -\f 9{640^2}\sum_{k=0}^{p-1}\f{\b{2k}k^2\b{4k}{2k}}{12^{4k}(k+1)}
-\f {223}{256000}\Ls p5p\mod{p^2}\ \t{if $\ls{-10}p=-1$.}
\endcases
\endalign$$
\endpro
\par{\bf Remark 2.18} Let $p>5$ be a prime. In [Su1], Z.W. Sun
conjectured
$$\sum_{k=0}^{p-1}(10k+1)\f{\b{2k}k^2\b{4k}{2k}}{12^{4k}}\e \Ls{-2}pp\mod{p^3},$$
and the congruence for
$\sum_{k=0}^{p-1}\f{\b{2k}k^2\b{4k}{2k}}{12^{4k}}$ modulo $p^2$. In
[S11], the author conjectured the congruence for
$\sum_{k=0}^{p-1}\f{\b{2k}k^2\b{4k}{2k}}{12^{4k}}$ modulo $p^3$.
Suppose that $\ls{-10}p=-1$ and
$$R_{40}(p)=\b{(p-1)/2}{[7p/40]}\b{(p-1)/2}{[9p/40]}
\b{[3p/40]}{[p/40]}\b{[19p/40]}{[p/20]}^{-1}.$$ In [S11], the author
conjectured that
$$\sum_{k=0}^{p-1}\f{\b{2k}k^2\b{4k}{2k}}{12^{4k}(k+1)}
\e\cases -\f{21}{5}R_{40}(p)\mod p&\t{if $p\e 3\mod {40}$,}
\\-\f{4446}{155}R_{40}(p)\mod p&\t{if $p\e 17\mod {40}$,}
\\-\f{189}{5}R_{40}(p)\mod p&\t{if $p\e 21\mod {40}$,}
\\-\f{702}{5}R_{40}(p)\mod p&\t{if $p\e 27\mod {40}$,}
\\\f{66}{5}R_{40}(p)\mod p&\t{if $p\e 29\mod {40}$,}
\\\f{1026}{5}R_{40}(p)\mod p&\t{if $p\e 31\mod {40}$,}
\\-\f{462}{5}R_{40}(p)\mod p&\t{if $p\e 33\mod {40}$,}
\\-\f{858}{85}R_{40}(p)\mod p&\t{if $p\e 39\mod {40}.$}
\endcases$$

\pro{Conjecture 2.35} Let $p>5$ be a prime. Then
$$\sum_{k=0}^{p-1}k^2\f{\b{2k}k^2\b{3k}k}{(-8640)^k}
\e\cases \f 4{729}x^2-\f{35}{1458}p+\f{p^2}{486x^2}\mod {p^3}
\\\qq\qq\qq\q\t{if $p\e 1,4\mod {15}$ and so $4p=x^2+75y^2$,}
\\-\f 4{243}x^2+\f{35}{1458}p-\f {p^2}{1458x^2}\mod {p^3}
\\\qq\qq\qq\q\t{if $p\e 7,13\mod {15}$ and so $4p=3x^2+25y^2$,}
\\\f 1{729}\sum_{k=0}^{p-1}\f{\b{2k}k^2\b{3k}k}{(-8640)^k(k+1)}
+\Ls p5\f p{1458}\mod {p^2}\q \t{if $3\mid p-2$}
\endcases$$
and
$$\sum_{k=0}^{p-1}k^3\f{\b{2k}k^2\b{3k}k}{(-8640)^k}
\e\cases -\f {19}{39366}x^2+\f{17}{78732}p+\f{p^2}{2916x^2}\mod
{p^3}
\\\qq\qq\qq\q\t{if $p\e 1,4\mod {15}$ and so $4p=x^2+75y^2$,}
\\\f {19}{13122}x^2-\f{17}{78732}p-\f {p^2}{8748x^2}\mod {p^3}
\\\qq\qq\qq\q\t{if $p\e 7,13\mod {15}$ and so $4p=3x^2+25y^2$,}
\\-\f 1{39366}\sum_{k=0}^{p-1}\f{\b{2k}k^2\b{3k}k}{(-8640)^k(k+1)}
+\Ls p5\f {77p}{78732}\mod {p^2}\q \t{if $3\mid p-2$.}
\endcases$$

\endpro
\par{\bf Remark 2.19} Let $p>5$ be a prime. In [S3] the author
conjectured that
$$\sum_{k=0}^{p-1}(9k+1)\f{\b{2k}k^2\b{3k}k}{(-8640)^k}\e
\Ls{-15}pp\mod {p^3}.$$ In [S11], the author conjectured that for
$p\e 1\mod 3$,
$$\sum_{k=0}^{p-1}\f{\b{2k}k^2\b{3k}k}{(-8640)^k}\e
\cases x^2-2p-\f{p^2}{x^2}\mod{p^3}&\t{if $p\e 1,4\mod {15}$ and so
$4p=x^2+75y^2$,}
\\-3x^2+2p+\f{p^2}{3x^2}\mod{p^3}&\t{if $p\e 7,13\mod {15}$ and so
$4p=3x^2+25y^2$.}
\endcases$$

\pro{Conjecture 2.36} Let $p$ be a prime with $p\not=2,3,17$. Then
$$\sum_{k=0}^{p-1}k^2\f{\b{2k}k^2\b{3k}k}{(-1728)^k}
\e\cases \f 2{289}x^2-\f{191}{5202}p+\f{47p^2}{5202x^2}\mod {p^3}
\\\qq\qq\qq\qq\t{if $\ls p3=\ls p{17}=1$ and so $4p=x^2+51y^2$,}
\\-\f 6{289}x^2+\f{191}{5202}p-\f {47p^2}{15606x^2}\mod {p^3}
\\\qq\qq\qq\qq\t{if $\ls p3=\ls p{17}=-1$ and so $4p=3x^2+17y^2$,}
\\\f 1{153}\sum_{k=0}^{p-1}\f{\b{2k}k^2\b{3k}k}{(-1728)^k(k+1)}
-\f 7{1734}\Ls p3p\mod {p^2}\  \t{if $\ls{-51}p=-1$}
\endcases$$
and
$$\sum_{k=0}^{p-1}k^3\f{\b{2k}k^2\b{3k}k}{(-1728)^k}
\e\cases -\f 1{88434}\big(\f
{349}3x^2+\f{209}{2}p+\f{49p^2}{6x^2}\big)\mod {p^3}
\\\qq\qq\qq\qq\qq\t{if $\ls p3=\ls p{17}=1$ and so $4p=x^2+51y^2$,}
\\\f 1{88434}\big(349x^2+\f{209}{2}p+\f {49p^2}{18x^2}\big)\mod {p^3}
\\\qq\qq\qq\qq\qq\t{if $\ls p3=\ls p{17}=-1$ and so $4p=3x^2+17y^2$,}
\\-\f 1{1734}\sum_{k=0}^{p-1}\f{\b{2k}k^2\b{3k}k}{(-1728)^k(k+1)}
-\f {2905}{530604}\Ls p3p\mod {p^2}\  \t{if $\ls{-51}p=-1$.}
\endcases$$
\endpro

\pro{Conjecture 2.37} Let $p$ be a prime with $p\not=2,3,41$. Then
$$\sum_{k=0}^{p-1}k^2\f{\b{2k}k^2\b{3k}k}{(-48)^{3k}}
\e\cases \f
1{630375}\big(2299x^2-\f{55837}{6}p+\f{661p^2}{6x^2}\big)\mod {p^3}
\\\qq\qq\qq\qq\t{if $\ls p3=\ls p{41}=1$ and so $4p=x^2+123y^2$,}
\\\f 1{630375}\big(-6897x^2+\f{55837}{6}p-\f {661p^2}{18x^2}\big)\mod {p^3}
\\\qq\qq\qq\qq\t{if $\ls p3=\ls p{41}=-1$ and so $4p=3x^2+41y^2$,}
\\\f 1{9225}\sum_{k=0}^{p-1}\f{\b{2k}k^2\b{3k}k}{(-48)^{3k}(k+1)}
-\f {53}{1260750}\Ls p3p\mod {p^2}\  \t{if $\ls{-123}p=-1$}
\endcases$$
and
$$\sum_{k=0}^{p-1}k^3\f{\b{2k}k^2\b{3k}k}{(-48)^{3k}}
\e\cases \f 1{6250\cdot
123^3}\big(-616598x^2+\f{1308183}{2}p+\f{1103101p^2}{2x^2}\big)\mod
{p^3}
\\\qq\qq\qq\qq\t{if $\ls p3=\ls p{41}=1$ and so $4p=x^2+123y^2$,}
\\\f 1{6250\cdot
123^3}\big(1849794x^2-\f{1308183}{2}p-\f {1103101p^2}{6x^2}\big)\mod
{p^3}
\\\qq\qq\qq\qq\t{if $\ls p3=\ls p{41}=-1$ and so $4p=3x^2+41y^2$,}
\\-\f 1{6303750}\big(\sum_{k=0}^{p-1}\f{\b{2k}k^2\b{3k}k}{(-48)^{3k}(k+1)}
+\f {1411019}{3690}\ls p3p\big)\mod {p^2}\  \t{if $\ls{-123}p=-1$.}
\endcases$$
\endpro

\par{\bf Remark 2.20} Let $p>3$ be a prime. In [Su4], Z.W. Sun
made conjectures equivalent to
$$\align&\sum_{k=0}^{p-1}(51k+7)\f{\b{2k}k^2\b{3k}k}{(-1728)^k}\e
7\Ls p3p+5p^3U_{p-3}\mod {p^4},
\\&\sum_{k=0}^{p-1}(615k+53)\f{\b{2k}k^2\b{3k}k}{(-48)^{3k}}
\e 53\Ls p3p+\f 52p^3U_{p-3}\mod {p^4}.
\endalign$$
For the congruences for
$\sum_{k=0}^{p-1}\f{\b{2k}k^2\b{3k}k}{(-1728)^k}$ and
$\sum_{k=0}^{p-1}\f{\b{2k}k^2\b{3k}k}{(-48)^{3k}}$ modulo $p^3$ see
[S11, Conjectures 5.39-5.40]. The corresponding congruences modulo
$p^2$ were conjectured by Z.W. Sun in [Su4]. \vskip0.2cm
\pro{Conjecture 2.38} Let $p$ be a prime with $p\not=2,3,5,13$. Then
$$\sum_{k=0}^{p-1}k^2\f{\b{2k}k^2\b{4k}{2k}}{(-82944)^{k}}
\e\cases \f
1{84500}\big(1271x^2-\f{5233}{4}p+\f{149p^2}{16x^2}\big)\mod {p^3}
\\\qq\qq\qq\q\t{if $\ls {-1}p=\ls {13}p=1$ and so $p=x^2+13y^2$,}
\\\f 1{84500}\big(-\f{1271}2x^2+\f{5233}{4}p-\f {149p^2}{8x^2}\big)\mod {p^3}
\\\qq\qq\qq\q\t{if $\ls {-1}p=\ls {13}p=-1$ and so $2p=x^2+13y^2$,}
\\\f 3{10400}\sum_{k=0}^{p-1}\f{\b{2k}k^2\b{4k}{2k}}{(-82944)^{k}(k+1)}
-\f {23}{169000}\Ls {-1}pp\mod {p^2}\\ \qq\qq\qq\q \t{if
$\ls{-13}p=-1$}
\endcases$$
and
$$\sum_{k=0}^{p-1}k^3\f{\b{2k}k^2\b{4k}{2k}}{(-82944)^{k}}
\e\cases \f
1{109850000}\big(-52151x^2+\f{38223}{4}p+\f{43631p^2}{16x^2}\big)\mod
{p^3}
\\\qq\qq\qq\q\t{if $\ls {-1}p=\ls {13}p=1$ and so $p=x^2+13y^2$,}
\\\f 1{109850000}\big(\f{52151}2x^2-\f{38223}{4}p-\f {43631p^2}{8x^2}\big)\mod {p^3}
\\\qq\qq\qq\q\t{if $\ls {-1}p=\ls {13}p=-1$ and so $2p=x^2+13y^2$,}
\\-\f 1{2600^2}\Big(9\sum_{k=0}^{p-1}\f{\b{2k}k^2\b{4k}{2k}}{(-82944)^{k}(k+1)}
+\f {81949}{65}\Ls {-1}pp\Big)\mod {p^2}\\ \qq\qq\qq\q \t{if
$\ls{-13}p=-1$.}
\endcases$$
\endpro
\par{\bf Remark 2.21} Let $p$ be a prime with $p\not=2,3,13$. Z.W.
Sun conjectured the congruence for
$\sum_{k=0}^{p-1}\f{\b{2k}k^2\b{4k}{2k}}{(-82944)^k}$ modulo $p^2$
and
$$\sum_{k=0}^{p-1}(260k+23)\f{\b{2k}k^2\b{4k}{2k}}{(-82944)^k}\e
23\Ls{-1}pp+\f 53p^3E_{p-3}\mod {p^4}.$$ See arXiv:0911.5665v59. For
the congruence for
$\sum_{k=0}^{p-1}\f{\b{2k}k^2\b{4k}{2k}}{(-82944)^k}$ modulo $p^3$
see [S11, Conjecture 5.42].

\section*{3. Conjectures for congruences involving Ap\'ery-like numbers}
\par\q
With the help of Maple, we discover the following conjectures
involving Ap\'ery-like numbers.
 \pro{Conjecture 3.1} Let $p$ be an odd prime. Then
$$
\align &\sum_{n=0}^{p-1}n^2A_n \e\cases
\f{15}{16}x^2-\f{31}{32}p+\f{p^2}{128x^2}\mod {p^3}&\t{if
$p=x^2+2y^2\e 1,3\mod 8$,}
\\ -\f 3{64}R_2(p)-\f p2\mod {p^2}&\t{if $p\e 5,7\mod 8$,}\endcases
\\&\sum_{n=0}^{p-1}n^3A_n \e\cases
-\f{13}{32}x^2+\f{37}{64}p-\f{19p^2}{256x^2}\mod {p^3}&\t{if
$p=x^2+2y^2\e 1,3\mod 8$,}
\\ \f 9{128}R_2(p)+\f 38p\mod {p^2}&\t{if $p\e 5,7\mod 8$.}\endcases
\endalign$$\endpro

\pro{Conjecture 3.2} Let $p$ be a prime with $p>3$. Then
$$
\align &\sum_{n=0}^{p-1}(-1)^nn^2A_n \e\cases
\f{8}{9}x^2-\f{17}{18}p+\f{p^2}{72x^2}\mod {p^3}&\t{if $p=x^2+3y^2\e
1\mod 3$,}
\\ \f 2{9}R_3(p)+\f p2\mod {p^2}&\t{if $p\e 2\mod 3$,}\endcases
\\&\sum_{n=0}^{p-1}(-1)^nn^3A_n \e\cases
-\f 13x^2+\f p2-\f{p^2}{12x^2}\mod {p^3}&\t{if $p=x^2+3y^2\e 1\mod
3$,}
\\ -\f 13R_3(p)-\f 13p\mod {p^2}&\t{if $p\e 2\mod 3$.}\endcases
\endalign$$\endpro
\par {\bf Remark 3.1} Let $p>3$ be a prime. In [Su2] Z.W. Sun showed
that
$$\sum_{n=0}^{p-1}(2n+1)A_n\e p+\f 76p^4B_{p-3}\mod {p^5}.$$
In [GuZ], Guo and Zeng proved Z.W. Sun's conjecture:
$$\sum_{n=0}^{p-1}(2n+1)(-1)^nA_n\e \Ls p3p\mod {p^3}.$$
For the conjectures on $\sum_{n=0}^{p-1}A_n$ and
$\sum_{n=0}^{p-1}(-1)^nA_n$ modulo $p^3$ see [S7, Conjecture 4.12].
The corresponding congruences modulo $p^2$ were conjectured by Z.W.
Sun earlier.\vskip0.2cm

\pro{Conjecture 3.3} Let $p$ be a prime with $p>3$. Then
$$
\align &\sum_{n=0}^{p-1}n^2\f{D_n}{(-2)^n} \\&\e\cases
\f{40}{27}x^2-\f{20+26(-1)^{(p-1)/2}}{27}p+\f{p^2}{18x^2}\mod
{p^3}&\t{if $p=x^2+3y^2\e 1\mod 3$,}
\\ \f 8{27}R_3(p)-\f{26}{27}(-1)^{\f{p-1}2}p\mod {p^2}&\t{if $p\e 2\mod 3$,}
\endcases
\\&\sum_{n=0}^{p-1}n^3\f{D_n}{(-2)^n}\\& \e\cases
-\f{40}{81}x^2+\f{20+68(-1)^{(p-1)/2}}{81}p-\f{17p^2}{54x^2}\mod
{p^3}&\t{if $p=x^2+3y^2\e 1\mod 3$,}
\\ -\f {56}{81}R_3(p)+\f{68}{81}(-1)^{\f{p-1}2}p\mod {p^2}&\t{if $p\e 2\mod 3$.}
\endcases
\endalign$$\endpro

\pro{Conjecture 3.4} Let $p$ be a prime with $p>3$. Then
$$
\align &\sum_{n=0}^{p-1}n^2\f{D_n}{(-32)^n} \\&\e\cases
\f{4}{27}x^2-\f{2+5(-1)^{(p-1)/2}}{27}p+\f{p^2}{36x^2}\mod
{p^3}&\t{if $p=x^2+3y^2\e 1\mod 3$,}
\\ \f 8{27}R_3(p)-\f{5}{27}(-1)^{\f{p-1}2}p\mod {p^2}&\t{if $p\e 2\mod 3$,}
\endcases
\\&\sum_{n=0}^{p-1}n^3\f{D_n}{(-32)^n}\\& \e\cases
\f{4}{81}x^2-\f{2+2(-1)^{(p-1)/2}}{81}p-\f{p^2}{36x^2}\mod
{p^3}&\t{if $p=x^2+3y^2\e 1\mod 3$,}
\\ -\f {16}{81}R_3(p)-\f{2}{81}(-1)^{\f{p-1}2}p\mod {p^2}&\t{if $p\e 2\mod 3$.}
\endcases
\endalign$$\endpro

\pro{Conjecture 3.5} Let $p$ be a prime with $p>3$. Then
$$
\align &\sum_{n=0}^{p-1}n^2\f{D_n}{4^n} \e\cases
\f{16}{9}x^2-\f{8}{9}p-\f{7p^2}{18x^2}\mod {p^3}&\t{if $p=x^2+3y^2\e
1\mod 3$,}
\\ -\f {20}9R_3(p)\mod {p^2}&\t{if $p\e 2\mod 3$,}
\endcases
\\&\sum_{n=0}^{p-1}n^3\f{D_n}{4^n} \e\cases
-\f{64}{45}x^2+\f{32}{45}p+\f{43p^2}{90x^2}\mod {p^3}&\t{if
$p=x^2+3y^2\e 1\mod 3$,}
\\ \f {28}{9}R_3(p)\mod {p^2}&\t{if $p\e 2\mod 3$ and $p\not=5$.}
\endcases
\endalign$$\endpro

\pro{Conjecture 3.6} Let $p$ be a prime with $p>3$. Then
$$
\align &\sum_{n=0}^{p-1}n^2\f{D_n}{16^n} \e\cases
\f{4}{9}x^2-\f{2}{9}p-\f{p^2}{18x^2}\mod {p^3}&\t{if $p=x^2+3y^2\e
1\mod 3$,}
\\ \f {4}9R_3(p)\mod {p^2}&\t{if $p\e 2\mod 3$,}
\endcases
\\&\sum_{n=0}^{p-1}n^3\f{D_n}{16^n} \e\cases
\f{4}{45}x^2-\f{2}{45}p+\f{p^2}{45x^2}\mod {p^3}&\t{if $p=x^2+3y^2\e
1\mod 3$,}
\\ -\f {4}{9}R_3(p)\mod {p^2}&\t{if $p\e 2\mod 3$.}
\endcases
\endalign$$\endpro

\pro{Conjecture 3.7} Let $p$ be an odd prime. Then
$$
\align &\sum_{n=0}^{p-1}n^2\f{D_n}{8^n} \e\cases
\f{3}{2}x^2-\f{5}{4}p-\f{p^2}{16x^2}\mod {p^3}&\t{if $p=x^2+2y^2\e
1,3\mod 8$,}
\\ -\f 3{8}R_2(p)-\f p2\mod {p^2}&\t{if $p\e 5,7\mod 8$,}\endcases
\\&\sum_{n=0}^{p-1}n^3\f{D_n}{8^n} \e\cases
-\f{5}{4}x^2+\f{17}{8}p+\f{p^2}{32x^2}\mod {p^3}&\t{if $p=x^2+2y^2\e
1,3\mod 8$,}
\\ \f 9{16}R_2(p)+\f 32p\mod {p^2}&\t{if $p\e 5,7\mod 8$.}\endcases
\endalign$$\endpro

\pro{Conjecture 3.8} Let $p$ be a prime with $p>5$. Then
$$
\align &\sum_{n=0}^{p-1}n^2D_n \\&\e\cases
\f{592}{225}x^2-\f{656}{225}p+\f{16p^2}{225x^2}\mod {p^3}&\t{if
$p=x^2+15y^2\e 1,19\mod {30}$,}
\\ -\f {592}{75}x^2+\f{656}{225}p-\f{16p^2}{675x^2}\mod {p^3}
&\t{if $p=3x^2+5y^2\e 17,23\mod {30}$,}
\\\f{16}{45}\sum_{k=0}^{p-1}\f{\b{2k}k^2\b{3k}k}{(-27)^k(k+1)}
-\f{296}{225}\ls p3p\mod {p^2}&\t{if $p\e 7,11,13,29\mod {30}$},
\endcases
\\&\sum_{n=0}^{p-1}n^3D_n \\&\e\cases
-\f{304}{125}x^2+\f{1456}{375}p-\f{87p^2}{125x^2}\mod {p^3}&\t{if
$p=x^2+15y^2\e 1,19\mod {30}$,}
\\ \f {912}{125}x^2-\f{1456}{375}p+\f{29p^2}{125x^2}\mod {p^3}
&\t{if $p=3x^2+5y^2\e 17,23\mod {30}$,}
\\-\f{32}{25}\sum_{k=0}^{p-1}\f{\b{2k}k^2\b{3k}k}{(-27)^k(k+1)}
+\f{616}{375}\ls p3p\mod {p^2}&\t{if $p\e 7,11,13,29\mod {30}$}.
\endcases
\endalign$$\endpro

\pro{Conjecture 3.9} Let $p$ be a prime with $p>5$. Then
$$
\align &\sum_{n=0}^{p-1}n^2\f{D_n}{64^n} \\&\e\cases
\f{52}{225}x^2-\f{26}{225}p-\f{13p^2}{450x^2}\mod {p^3}&\t{if
$p=x^2+15y^2\e 1,19\mod {30}$,}
\\ -\f {52}{75}x^2+\f{26}{225}p+\f{13p^2}{1350x^2}\mod {p^3}
&\t{if $p=3x^2+5y^2\e 17,23\mod {30}$,}
\\\f{16}{45}\sum_{k=0}^{p-1}\f{\b{2k}k^2\b{3k}k}{(-27)^k(k+1)}
+\f{64}{225}\ls p3p\mod {p^2}&\t{if $p\e 7,11,13,29\mod {30}$},
\endcases
\\&\sum_{n=0}^{p-1}n^3\f{D_n}{64^n} \\&\e\cases
\f{52}{375}x^2-\f{p}{375}-\f{13p^2}{750x^2}\mod {p^3}&\t{if
$p=x^2+15y^2\e 1,19\mod {30}$,}
\\ -\f {52}{125}x^2+\f{p}{375}+\f{13p^2}{2250x^2}\mod {p^3}
&\t{if $p=3x^2+5y^2\e 17,23\mod {30}$,}
\\\f{16}{75}\sum_{k=0}^{p-1}\f{\b{2k}k^2\b{3k}k}{(-27)^k(k+1)}
+\f{89}{375}\ls p3p\mod {p^2}&\t{if $p\e 7,11,13,29\mod {30}$}.
\endcases
\endalign$$\endpro

\pro{Conjecture 3.10} Let $p>3$ be a prime. Then
$$\align &\sum_{n=0}^{p-1}n^2\f{D_n}{(-8)^n} \\&\e\cases \f{11}{18}x^2
-\f{29}{36}p+\f{7}{144x^2}p^2\mod {p^3} &\t{if $p=x^2+6y^2\e 1,7\mod
{24},$}
\\\f{11}{9}x^2
+\f{7}{36}p+\f{7}{288x^2}p^2\mod {p^3} &\t{if $p=2x^2+3y^2\e
5,11\mod {24},$}
\\-\f 19\sum_{k=0}^{p-1}\f{\b{2k}k^2\b{3k}k}{216^k(k+1)}-\f{17}{36}
\Ls p3p\mod {p^2}&\t{if $p\e 13,17,19,23\mod {24}$}\endcases
\endalign$$
and
$$\align &\sum_{n=0}^{p-1}n^3\f{D_n}{(-8)^n} \\&\e\cases \f{x^2}{12}
+\f p8-\f{13p^2}{96x^2}\mod {p^3} &\t{if $p=x^2+6y^2\e 1,7\mod
{24},$}
\\\f{x^2}{6}
-\f{5p}{24}-\f{13p^2}{192x^2}\mod {p^3} &\t{if $p=2x^2+3y^2\e
5,11\mod {24},$}
\\\f 16\sum_{k=0}^{p-1}\f{\b{2k}k^2\b{3k}k}{216^k(k+1)}+\f{1}{8}
\Ls p3p\mod {p^2}&\t{if $p\e 13,17,19,23\mod {24}$.}\endcases
\endalign$$
\endpro
\par{\bf Remark 3.2} Let $p>3$ be a prime. In [S7], the author conjectured the congruences modulo $p^3$ for
$\sum_{n=0}^{p-1}\f{D_n}{m^n}$ $(m\in\{1,-2,4,8,-8,16,-32,64\})$.
The corresponding congruences modulo $p^2$ were conjectured by Z.W.
Sun in [Su4]. In [Su4], Z.W. Sun also made conjectures for
$$\sum_{n=0}^{p-1}(5n+4)D_n,\q
\sum_{n=0}^{p-1}(3n+2)\f{D_n}{(-2)^n},\q
\sum_{n=0}^{p-1}(2n+1)\f{D_n}{(-8)^n}$$ modulo $p^3$ and
 $\sum_{n=0}^{p-1}(2n+1)\f{D_n}{8^n}$ modulo $p^4$.
 In [S6], the author conjectured that
$$\align &\sum_{n=0}^{p-1}(5n+4)D_n\e 4p\Ls p3+28p^3U_{p-3}\mod
{p^4},
\\&\sum_{n=0}^{p-1}(3n+2)\f{D_n}{(-2)^n}\e
2(-1)^{\f{p-1}2}p+6p^3E_{p-3}\mod{p^4},
\\&\sum_{n=0}^{p-1}(2n+1)\f{D_n}{(-8)^n}\e
p\Ls p3+\f 52p^3U_{p-3}\mod{p^4},
\\&\sum_{n=0}^{p-1}(2n+1)\f{D_n}{8^n}\e
p+\f {35}{24}p^4B_{p-3}\mod{p^5},
\\&\sum_{n=0}^{p-1}(5n+1)\f{D_n}{64^n}\e
p\Ls p3-2p^3U_{p-3}\mod{p^4}.
\endalign$$
In [Liu3], Liu proved the congruences for
$\sum_{n=0}^{p-1}(3n+2)\f{D_n}{(-2)^n}$ and
$\sum_{n=0}^{p-1}(3n+1)\f{D_n}{(-32)^n}$ modulo $p^4$. Few days ago
Mao and Liu [ML] proved the author's conjecture
$$\align&\sum_{n=0}^{p-1}\f{D_n}{4^n}\e \f{3\sls p3-1}2
\sum_{n=0}^{p-1}\f{D_n}{16^n}
\\&\q\e\cases 4x^2-2p-\f{p^2}{4x^2}\mod {p^3}
&\t{if $p=x^2+3y^2\e 1\mod 3$,}
\\\f{p^2}2\b{(p-1)/2}{(p-5)/6}^{-2}\mod {p^3}&\t{if $p\e 2\mod 3$.}
\endcases\endalign$$
The corresponding congruences modulo $p^2$ were conjectured by Z.W.
Sun in [Su4] and proved by the author in [S10].

 \vskip0.2cm \pro{Conjecture 3.11} Let $p>3$ be a
prime. Then
$$
\align &\sum_{n=0}^{p-1}n^2b_n \\&\e\cases
\f{33}{16}x^2-\f{33+40\sls p3}{32}p+\f{7p^2}{128x^2}\mod {p^3}&\t{if
$p=x^2+2y^2\e 1,3\mod 8$,}
\\ \f 3{64}R_2(p)-\f 54\ls p3p\mod {p^2}&\t{if $p\e 5,7\mod 8$,}\endcases
\\&\sum_{n=0}^{p-1}n^3b_n \\&\e\cases
-\f{75}{64}x^2+\f{75+184\sls p3}{128}p-\f{221p^2}{512x^2}\mod
{p^3}&\t{if $p=x^2+2y^2\e 1,3\mod 8$,}
\\ -\f {33}{256}R_2(p)+\f {23}{16}\ls p3p\mod {p^2}&\t{if $p\e 5,7\mod 8$.}\endcases
\endalign$$\endpro

\pro{Conjecture 3.12} Let $p>3$ be a prime. Then
$$
\align &\sum_{n=0}^{p-1}n^2\f{b_n}{81^n} \\&\e\cases
\f{x^2}{16}-\f{3+8\sls p3}{96}p+\f{5p^2}{384x^2}\mod {p^3}&\t{if
$p=x^2+2y^2\e 1,3\mod 8$,}
\\ \f 3{64}R_2(p)-\f 1{12}\ls p3p\mod {p^2}&\t{if $p\e 5,7\mod 8$,}\endcases
\\&\sum_{n=0}^{p-1}n^3\f{b_n}{81^n} \\&\e\cases
-\f{1}{64}x^2+\f{3-8\sls p3}{384}p-\f{5p^2}{1536x^2}\mod {p^3}&\t{if
$p=x^2+2y^2\e 1,3\mod 8$,}
\\ -\f {3}{256}R_2(p)-\f {1}{48}\ls p3p\mod {p^2}&\t{if $p\e 5,7\mod 8$.}\endcases
\endalign$$\endpro

\pro{Conjecture 3.13} Let $p$ be a prime with $p>3$. Then
$$
\align &\sum_{n=0}^{p-1}n^2\f{b_n}{(-9)^n} \e\cases -\f
p2+\f{p^2}{8x^2}\mod {p^3}&\t{if $p=x^2+3y^2\e 1\mod 3$,}
\\ 2R_3(p)+\f p2\mod {p^2}&\t{if $p\e 2\mod 3$,}
\endcases
\\&\sum_{n=0}^{p-1}n^3\f{b_n}{(-9)^n} \e\cases
x^2-\f 32p-\f{p^2}{4x^2}\mod {p^3}&\t{if $p=x^2+3y^2\e 1\mod 3$,}
\\ -3R_3(p)+p\mod {p^2}&\t{if $p\e 2\mod 3$ and $p\not=5$.}
\endcases
\endalign$$\endpro

\pro{Conjecture 3.14} Let $p$ be a prime with $p>3$. Then
$$
\align &\sum_{n=0}^{p-1}n^2\f{b_n}{(-3)^n} \\&\e\cases \f
94x^2-\f{21}8p+\f{3p^2}{32x^2}\mod {p^3}&\t{if $12\mid p-1$ and so
$p=x^2+9y^2$,}
\\-\f 98x^2+\f{21}8p-\f{3p^2}{16x^2}\mod {p^3}&\t{if $12\mid p-5$ and so
$2p=x^2+9y^2$,}
\\\f
3{128}\sum_{k=0}^{p-1}\f{\b{2k}k^2\b{4k}{2k}}{(-12288)^k(k+1)}-\f{87}{64}
\ls p3p\mod {p^2}&\t{if $p\e 3\mod 4$,}
\endcases
\\&\sum_{n=0}^{p-1}n^3\f{b_n}{(-3)^n} \\&\e\cases -\f
{27}{16}x^2+\f{105}{32}p-\f{99p^2}{128x^2}\mod {p^3}&\t{if $12\mid
p-1$ and so $p=x^2+9y^2$,}
\\\f {27}{32}x^2-\f{105}{32}p+\f{99p^2}{64x^2}\mod {p^3}&\t{if $12\mid p-5$ and so
$2p=x^2+9y^2$,}
\\-\f {45}{512}\sum_{k=0}^{p-1}\f{\b{2k}k^2\b{4k}{2k}}{(-12288)^k(k+1)}+
\f{489}{256} \ls p3p\mod {p^2}&\t{if $p\e 3\mod 4$.}
\endcases
\endalign$$\endpro

\pro{Conjecture 3.15} Let $p$ be a prime with $p>3$. Then
$$
\align &\sum_{n=0}^{p-1}n^2\f{b_n}{(-27)^n} \\&\e\cases \f
{1}{4}x^2-\f{1}{8}p-\f{p^2}{32x^2}\mod {p^3}&\t{if $12\mid p-1$ and
so $p=x^2+9y^2$,}
\\-\f 18x^2+\f 18p+\f{p^2}{16x^2}\mod {p^3}&\t{if $12\mid p-5$ and so
$2p=x^2+9y^2$,}
\\\f {3}{128}\sum_{k=0}^{p-1}\f{\b{2k}k^2\b{4k}{2k}}{(-12288)^k(k+1)}+
\f{9}{64} \ls p3p\mod {p^2}&\t{if $p\e 3\mod 4$,}
\endcases
\\&\sum_{n=0}^{p-1}n^3\f{b_n}{(-27)^n} \\&\e\cases -\f
{x^2}{16}+\f{3}{32}p-\f{p^2}{128x^2}\mod {p^3}&\t{if $12\mid p-1$
and so $p=x^2+9y^2$,}
\\\f{x^2}{32}-\f{3}{32}p+\f{p^2}{64x^2}\mod {p^3}&\t{if $12\mid p-5$ and so
$2p=x^2+9y^2$,}
\\\f
9{512}\sum_{k=0}^{p-1}\f{\b{2k}k^2\b{4k}{2k}}{(-12288)^k(k+1)}+\f{43}{256}
\ls p3p\mod {p^2}&\t{if $p\e 3\mod 4$.}
\endcases
\endalign$$\endpro

\pro{Conjecture 3.16} Let $p$ be a prime with $p>3$. Then
$$\align &\sum_{n=0}^{p-1}n^2\f{b_n}{9^n}
\\&\e\cases \f 9{16}x^2-\f{25}{32}p+\f{7p^2}{128x^2}
\mod {p^3}&\t{if $p=x^2+6y^2\e 1,7\mod {24}$,}
\\-\f 98x^2+\f{25}{32}p-\f{7p^2}{256x^2}\mod {p^3}&\t{if
$p=2x^2+3y^2\e 5,11\mod {24}$,}
\\\f 18\Ls p3\sum_{k=0}^{p-1}\f{\b{2k}k^2\b{3k}k}{216^k(k+1)}-\f
{1+16\sls p3}{32}p\mod {p^2}&\t{if $p\e 11,13,17,19\mod{24}$}
\endcases\endalign$$
and
$$\align &\sum_{n=0}^{p-1}n^3\f{b_n}{9^n}
\\&\e\cases \f 5{32}x^2+\f{3}{64}p-\f{37p^2}{256x^2}
\mod {p^3}&\t{if $p=x^2+6y^2\e 1,7\mod {24}$,}
\\-\f 5{16}x^2-\f{3}{64}p+\f{37p^2}{512x^2}\mod {p^3}&\t{if
$p=2x^2+3y^2\e 5,11\mod {24}$,}
\\-\f 3{16}\Ls p3\sum_{k=0}^{p-1}\f{\b{2k}k^2\b{3k}k}{216^k(k+1)}+\f
{3+8\sls p3}{64}p\mod {p^2}&\t{if $p\e 11,13,17,19\mod{24}$.}
\endcases\endalign$$
\endpro

\par{\bf Remark 3.3} Let $p>3$ be a prime. In [S7], the author
conjectured the congruences modulo $p^3$ for
$\sum_{n=0}^{p-1}\f{b_n}{m^n}$ in the cases $m=1,-3,9,-9,-27,81$. In
[S5], the author conjectured that
$$\align&\sum_{k=0}^{p-1}(4k+1)\f{b_k}{(-27)^k}
\e \sum_{k=0}^{p-1}(4k+1)\f{b_k}{81^k}\e
\sum_{k=0}^{p-1}(2k+1)\f{b_k}{(-9)^k} \e
\sum_{k=0}^{p-1}(2k+1)\f{b_k}{9^k}
\\&\e \f 13\sum_{k=0}^{p-1}(4k+3)b_k
\e \f 13\sum_{k=0}^{p-1}(4k+3)\f{b_k}{(-3)^k} \e \Ls p3p\mod{p^3}.
\endalign$$
The congruences for $\sum_{k=0}^{p-1}(4k+1)\f{b_k}{81^k}$,
$\sum_{k=0}^{p-1}(4k+3)b_k$ and
$\sum_{k=0}^{p-1}(2k+1)\f{b_k}{(-9)^k}$ modulo $p^3$ have been
confirmed by J.-C. Liu in [Liu2].

\par\q
\pro{Conjecture 3.17} Let $p$ be a prime with $p\not=2,7$. Then
$$
\align &\sum_{n=0}^{p-1}n^2T_n \e\cases
\f{80}{49}x^2-\f{40}{49}p-\f{10p^2}{49x^2}\mod {p^3}&\t{if
$p=x^2+7y^2\e 1,2,4\mod 7$,}
\\\f{16}7\sum_{k=0}^{(p-1)/2}\f{\b{2k}k^3}{k+1}+\f{128}{49}p\mod {p^2}
&\t{if $p\e 3,5,6\mod 7$,}
\endcases
\\&\sum_{n=0}^{p-1}n^3T_n \e\cases
\f{176}{343}x^2+\f{696}{343}p-\f{71p^2}{343x^2}\mod {p^3}&\t{if
$p=x^2+7y^2\e 1,2,4\mod 7$,}
\\\f{192}{49}\sum_{k=0}^{(p-1)/2}\f{\b{2k}k^3}{k+1}+\f{2320}{343}p\mod {p^2}
&\t{if $p\e 3,5,6\mod 7$.}
\endcases
\endalign$$
\endpro

\pro{Conjecture 3.18} Let $p$ be a prime with $p\not=2,7$. Then
$$
\align &\sum_{n=0}^{p-1}n^2\f{T_n}{16^n} \e\cases
\f{52}{49}x^2-\f{68}{49}p+\f{4p^2}{49x^2}\mod {p^3}&\t{if
$p=x^2+7y^2\e 1,2,4\mod 7$,}
\\\f{16}7\sum_{k=0}^{(p-1)/2}\f{\b{2k}k^3}{k+1}+\f{86}{49}p\mod {p^2}
&\t{if $p\e 3,5,6\mod 7$,}
\endcases
\\&\sum_{n=0}^{p-1}n^3\f{T_n}{16^n} \e\cases
-\f{876}{343}x^2+\f{1467}{343}p-\f{369p^2}{686x^2}\mod {p^3}&\t{if
$p=x^2+7y^2\e 1,2,4\mod 7$,}
\\-\f{528}{49}\sum_{k=0}^{(p-1)/2}\f{\b{2k}k^3}{k+1}-\f{3195}{343}p\mod {p^2}
&\t{if $p\e 3,5,6\mod 7$.}
\endcases
\endalign$$
\endpro

\pro{Conjecture 3.19} Let $p$ be an odd prime. Then
$$
\align &\sum_{n=0}^{p-1}n^2\f{T_n}{4^n} \e\cases -4y^2\mod
{p^3}&\t{if $p=x^2+4y^2\e 1\mod 4$,}
\\-\f 14R_1(p)-\f 12p\mod {p^2}
&\t{if $p\e 3\mod 4$,}
\endcases
\\&\sum_{n=0}^{p-1}n^3\f{T_n}{4^n} \e\cases
2y^2+\f p4+\f{p^2}{64y^2}\mod {p^3}&\t{if $p=x^2+4y^2\e 1\mod 4$,}
\\\f{3}{8}R_1(p)+\f 12p\mod {p^2}
&\t{if $p\e 3\mod 4$.}
\endcases
\endalign$$
\endpro

\pro{Conjecture 3.20} Let $p$ be an odd prime. Then
$$\align &\sum_{n=0}^{p-1}n^2\f{T_n}{(-4)^n}\\& \e\cases \f
34x^2-\f{3+4(-1)^{\f{p-1}2}}{8}p+\f{p^2}{32x^2}\mod {p^3}&\t{if
$p=x^2+2y^2\e 1,3\mod 8$,}
\\\f 1{16}R_2(p)-\f 12(-1)^{\f{p-1}2}p\mod {p^2}
&\t{if $p\e 5,7\mod 8$,}
\endcases
\\&\sum_{n=0}^{p-1}n^3\f{T_n}{(-4)^n}\\& \e\cases -\f
18x^2+\f{1+4(-1)^{\f{p-1}2}}{16}p-\f{7p^2}{64x^2}\mod {p^3}&\t{if
$p=x^2+2y^2\e 1,3\mod 8$,}
\\-\f 3{32}R_2(p)+\f 14(-1)^{\f{p-1}2}p\mod {p^2}
&\t{if $p\e 5,7\mod 8$.}
\endcases
\endalign$$\endpro

\par{\bf Remark 3.4} Let $p>3$ be a prime. In [S6], the author
proved that
$$\sum_{n=0}^{p-1}\f{T_n}{4^n} \e \cases
4x^2-2p-\f{p^2}{4x^2}\mod {p^3}&\t{if $p=x^2+4y^2\e 1\mod 4$,}
\\-\f{p^2}4\b{(p-3)/2}{(p-3)/4}^{-2}\mod {p^3}
&\t{if $p\e 3\mod 4$,}
\endcases$$
and made conjectures on congruences modulo $p^3$ for
$\sum_{n=0}^{p-1}\f{T_n}{m^n}$ in the cases $m=1,-4,16$. In [S6],
the author also proved that
$$\align&\sum_{n=0}^{p-1}(2n+1)\f{T_n}{4^n}\e p\mod {p^4},
\\&\sum_{n=0}^{p-1}(2n+1)\f{T_n}{(-4)^n}\e (-1)^{\f{p-1}2}p\mod
{p^3}\endalign$$ and conjectured that
$$\align
&\sum_{n=0}^{p-1}(7n+4)T_n\e 4p+\f{25}3p^4B_{p-3}\mod {p^5},
\\&\sum_{n=0}^{p-1}(7n+3)\f{T_n}{16^n}\e 3p+\f{25}{12}p^4B_{p-3}\mod {p^5}.
 \endalign$$

\par\q
\pro{Conjecture 3.21} Let $p$ be an odd prime. Then
$$
\align &\sum_{n=0}^{p-1}n^2\f{V_n}{8^n} \e\cases -24y^2\mod
{p^3}&\t{if $p=x^2+4y^2\e 1\mod 4$,}
\\\f 12R_1(p)+3p\mod {p^2}
&\t{if $p\e 3\mod 4$,}
\endcases
\\&\sum_{n=0}^{p-1}n^3\f{V_n}{8^n} \e\cases
-16x^2+18p-\f{5p^2}{4x^2}\mod {p^3}&\t{if $p=x^2+4y^2\e 1\mod 4$,}
\\-3R_1(p)-10p\mod {p^2}
&\t{if $p\e 3\mod 4$.}
\endcases
\endalign$$
\endpro

\pro{Conjecture 3.22} Let $p$ be an odd prime. Then
$$
\align &\sum_{n=0}^{p-1}n^2\f{V_n}{(-16)^n} \e\cases \f 12x^2-\f
34p+\f{p^2}{16x^2}\mod {p^3}&\t{if $p=x^2+4y^2\e 1\mod 4$,}
\\\f 18R_1(p)+\f 12p\mod {p^2}
&\t{if $p\e 3\mod 4$,}
\endcases
\\&\sum_{n=0}^{p-1}n^3\f{V_n}{(-16)^n} \e\cases
\f 14x^2-\f{5p^2}{32x^2}\mod {p^3}&\t{if $p=x^2+4y^2\e 1\mod 4$,}
\\-\f 3{16}R_1(p)-\f 18p\mod {p^2}
&\t{if $p\e 3\mod 4$.}
\endcases
\endalign$$
\endpro

\pro{Conjecture 3.23} Let $p$ be an odd prime. Then
$$
\align &\sum_{n=0}^{p-1}n^2\f{V_n}{32^n} \e\cases
2x^2-p-\f{p^2}{4x^2}\mod {p^3}&\t{if $p=x^2+4y^2\e 1\mod 4$,}
\\\f 12R_1(p)\mod {p^2}
&\t{if $p\e 3\mod 4$,}
\endcases
\\&\sum_{n=0}^{p-1}n^3\f{V_n}{32^n} \e\cases
6x^2-3p-\f{3p^2}{4x^2}\mod {p^3}&\t{if $p=x^2+4y^2\e 1\mod 4$,}
\\\f 3{2}R_1(p)\mod {p^2}
&\t{if $p\e 3\mod 4$.}
\endcases
\endalign$$
\endpro
\par{\bf Remark 3.5} Let $p$ be an odd prime. In [S8], the author
conjectured the congruence for $\sum_{n=0}^{p-1}\f{V_n}{m^n}\mod
{p^3}$ in the cases $m=8,-16,32$.\vskip0.2cm

\pro{Conjecture 3.24} Let $p>3$ be a prime. Then
$$\align &\sum_{n=0}^{p-1}n^2\f{V_n^{(3)}}{3^n}
\\&\e\cases \f
{13}{16}x^2-\f {27p}{8}+\f{p^2}{8x^2}\mod{p^3}&\t{if $3\mid p-1$ and
so $4p=x^2+27y^2$,}
\\\f 18(2p+1)\b{[2p/3]}{[p/3]}^2+\f 74p\mod {p^2}&\t{if $p\e 2\mod{3}$}
\endcases\endalign$$
and
$$\align &\sum_{n=0}^{p-1}n^3\f{V_n^{(3)}}{3^n}
\\&\e\cases -\f
{111}{128}x^2+\f {293p}{64}-\f{147p^2}{64x^2}\mod{p^3}&\t{if $3\mid
p-1$ and so $4p=x^2+27y^2$,}
\\-\f {27}{64}(2p+1)\b{[2p/3]}{[p/3]}^2-\f {91}{32}
p\mod {p^2}&\t{if $p\e 2\mod{3}$.}
\endcases\endalign$$
Moreover,
$$\sum_{n=0}^{p-1}(8n+7)\f{V_n^{(3)}}{3^n}
\e 7\Ls p3p+\f{278}3p^3U_{p-3}\mod {p^4}.$$
\endpro

\pro{Conjecture 3.25} Let $p>3$ be a prime. Then

$$\align &\sum_{n=0}^{p-1}n^2\f{V_n^{(3)}}{243^n}
\e\cases \f {1}{16}x^2-\f {1}{8}p-\f{p^2}{8x^2}\mod{p^3}&\t{if
$3\mid p-1$ and so $4p=x^2+27y^2$,}
\\\f 18(2p+1)\b{[2p/3]}{[p/3]}^2\mod {p^2}&\t{if $p\e 2\mod{3}$}
\endcases\endalign$$
and
$$\align &\sum_{n=0}^{p-1}n^3\f{V_n^{(3)}}{243^n}
\e\cases \f {7}{128}x^2-\f {5}{64}p-\f{5p^2}{64x^2}\mod{p^3}&\t{if
$3\mid p-1$ and so $4p=x^2+27y^2$,}
\\\f 3{64}(2p+1)\b{[2p/3]}{[p/3]}^2-\f p{32}\mod {p^2}&\t{if $p\e 2\mod{3}$.}
\endcases\endalign$$
Moreover,
$$\sum_{n=0}^{p-1}(8n+1)\f{V_n^{(3)}}{243^n}
\e \Ls p3p-\f{22}3p^3U_{p-3}\mod {p^4}.$$

\endpro

\pro{Conjecture 3.26} Let $p>3$ be a prime. Then
$$\sum_{n=0}^{p-1}n^2\f{V_n^{(3)}}{(-27)^n}
\e\cases \f {4}{9}x^2-\f {13}{18}p+\f{5p^2}{72x^2}\mod{p^3}&\t{if
$p=x^2+3y^2\e 1\mod 3$,}
\\\f 29R_3(p)+\f p2\mod {p^2}&\t{if $p\e 2\mod 3$}\endcases$$
and
$$\sum_{n=0}^{p-1}n^3\f{V_n^{(3)}}{(-27)^n}
\e\cases \f {1}{3}x^2-\f {p}{18}-\f{p^2}{6x^2}\mod{p^3}&\t{if
$p=x^2+3y^2\e 1\mod 3$,}
\\-\f 13R_3(p)-\f p9\mod {p^2}&\t{if $p\e 2\mod 3$.}\endcases$$

\endpro
\pro{Conjecture 3.27} Let $p$ be a prime with $p\not=2,3,7$. Then
$$\sum_{n=0}^{p-1}n^2V_n^{(4)}
\e\cases \f{4192}{1323}x^2-\f{4336}{1323}p+\f{4p^2}{147x^2}\mod
{p^3}&\t{if $p=x^2+7y^2\e 1,2,4\mod 7$,}
\\-\f 8{63}R_7(p)+\f{2048}{1323}p\mod{p^2}&\t{if $p\e 3,5,6\mod
7$}\endcases$$ and
$$\sum_{n=0}^{p-1}n^3V_n^{(4)}
\e\cases
-\f{238816}{83349}x^2+\f{107600}{27783}p-\f{42290p^2}{83349x^2}\mod
{p^3}&\t{if $p=x^2+7y^2\e 1,2,4\mod 7$,}
\\\f {512}{1323}R_7(p)-\f{166528}{83349}p\mod{p^2}&\t{if $p\e 3,5,6\mod
7$}.\endcases$$ Moreover,
$$\sum_{n=0}^{p-1}(9n+8)V_n^{(4)}\e 8\Ls p7p\mod {p^3}.$$
\endpro

\pro{Conjecture 3.28} Let $p$ be a prime with $p\not=2,3,7$. Then
$$\sum_{n=0}^{p-1}n^2\f{V_n^{(4)}}{4096^n}
\e\cases \f{76}{1323}x^2-\f{52}{1323}p-\f{2p^2}{441x^2}\mod
{p^3}&\t{if $p=x^2+7y^2\e 1,2,4\mod 7$,}
\\-\f 8{63}R_7(p)-\f{178}{1323}p\mod{p^2}&\t{if $p\e 3,5,6\mod
7$}\endcases$$ and
$$\sum_{n=0}^{p-1}n^3\f{V_n^{(4)}}{4096^n}
\e\cases
\f{2188}{83349}x^2-\f{269}{27783}p-\f{281p^2}{166698x^2}\mod
{p^3}&\t{if $p=x^2+7y^2\e 1,2,4\mod 7$,}
\\-\f {8}{1323}R_7(p)-\f{863}{83349}p\mod{p^2}&\t{if $p\e 3,5,6\mod
7$}.\endcases$$ Moreover,
$$\sum_{n=0}^{p-1}(9n+1)\f{V_n^{(4)}}{4096^n}\e \Ls p7p\mod {p^3}.$$
\endpro

\pro{Conjecture 3.29} Let $p$ be a prime with $p>3$. Then
$$\sum_{n=0}^{p-1}n^2\f{V_n^{(4)}}{16^n}
\e\cases \f{44}{9}x^2-\f{43}{9}p-\f{p^2}{36x^2}\mod {p^3}&\t{if
$p=x^2+3y^2\e 1\mod 3$,}
\\\f 2{9}R_3(p)+\f 73p\mod{p^2}&\t{if $p\e 2\mod
3$}\endcases$$ and
$$\sum_{n=0}^{p-1}n^3\f{V_n^{(4)}}{16^n}
\e\cases -\f{74}{9}x^2+\f{82}{9}p-\f{43p^2}{72x^2}\mod {p^3}&\t{if
$p=x^2+3y^2\e 1\mod 3$,}
\\-\f 8{9}R_3(p)-5p\mod{p^2}&\t{if $p\e 2\mod
3$}.\endcases$$ Moreover,
$$\sum_{n=0}^{p-1}(n+1)\f{V_n^{(4)}}{16^n}\e \Ls p3p+\f{35}2p^3U_{p-3}
\mod {p^4}.$$
\endpro

\pro{Conjecture 3.30} Let $p$ be a prime with $p>3$. Then
$$\sum_{n=0}^{p-1}n^2\f{V_n^{(4)}}{256^n}
\e\cases \f{8}{9}x^2-\f{4}{9}p-\f{p^2}{9x^2}\mod {p^3}&\t{if
$p=x^2+3y^2\e 1\mod 3$,}
\\\f 2{9}R_3(p)\mod{p^2}&\t{if $p\e 2\mod
3$}\endcases$$ and
$$\sum_{n=0}^{p-1}n^3\f{V_n^{(4)}}{256^n}
\e\cases \f{14}{9}x^2-\f{7}{9}p-\f{11p^2}{72x^2}\mod {p^3}&\t{if
$p=x^2+3y^2\e 1\mod 3$,}
\\\f 2{9}R_3(p)\mod{p^2}&\t{if $p\e 2\mod
3$}.\endcases$$ Moreover,
$$\sum_{n=0}^{p-1}\f{nV_n^{(4)}}{256^n}\e -\f{15}4p^3U_{p-3}\mod {p^4}.$$
\endpro

\pro{Conjecture 3.31} Let $p$ be a prime with $p>3$. Then
$$\sum_{n=0}^{p-1}n^2\f{V_n^{(4)}}{(-8)^n}
\e\cases \f{58}{27}x^2-\f{64}{27}p+\f{p^2}{18x^2}\mod {p^3}&\t{if
$p=x^2+4y^2\e 1\mod 4$,}
\\\f 1{18}R_1(p)+\f{35}{27}p\mod{p^2}&\t{if $p\e 3\mod
4$}\endcases$$ and
$$\sum_{n=0}^{p-1}n^3\f{V_n^{(4)}}{(-8)^n}
\e\cases -\f{260}{243}x^2+\f{160}{81}p-\f{427p^2}{972x^2}\mod
{p^3}&\t{if $p=x^2+4y^2\e 1\mod 4$,}
\\-\f 4{27}R_1(p)-\f{350}{243}p\mod{p^2}&\t{if $p\e 3\mod
4$}.\endcases$$ Moreover,
$$\sum_{n=0}^{p-1}(9n+7)\f{V_n^{(4)}}{(-8)^n}\e 7(-1)^{\f{p-1}2}p+79p^3
E_{p-3}\mod {p^4}.$$
\endpro

\pro{Conjecture 3.32} Let $p$ be a prime with $p>3$. Then
$$\sum_{n=0}^{p-1}n^2\f{V_n^{(4)}}{(-512)^n}
\e\cases -\f{2}{27}x^2-\f{p}{27}+\f{p^2}{36x^2}\mod {p^3}&\t{if
$p=x^2+4y^2\e 1\mod 4$,}
\\\f 1{18}R_1(p)+\f{2}{27}p\mod{p^2}&\t{if $p\e 3\mod
4$}\endcases$$ and
$$\sum_{n=0}^{p-1}n^3\f{V_n^{(4)}}{(-512)^n}
\e\cases -\f{10}{243}x^2-\f{p}{81}-\f{5p^2}{972x^2}\mod {p^3}&\t{if
$p=x^2+4y^2\e 1\mod 4$,}
\\-\f 1{54}R_1(p)+\f{8}{243}p\mod{p^2}&\t{if $p\e 3\mod
4$}.\endcases$$ Moreover,
$$\sum_{n=0}^{p-1}(9n+2)\f{V_n^{(4)}}{(-512)^n}\e 2(-1)^{\f{p-1}2}p+8p^3
E_{p-3}\mod {p^4}.$$
\endpro

\pro{Conjecture 3.33} Let $p$ be a prime with $p>3$. Then
$$\sum_{n=0}^{p-1}n^2\f{V_n^{(4)}}{(-64)^n}
\e\cases \f{3}{8}x^2-\f{11}{16}p+\f{5p^2}{64x^2}\mod {p^3}&\t{if
$p=x^2+2y^2\e 1,3\mod 8$,}
\\\f 1{32}R_2(p)+\f{p}{2}\mod{p^2}&\t{if $p\e 5,7\mod
8$}\endcases$$ and
$$\sum_{n=0}^{p-1}n^3\f{V_n^{(4)}}{(-64)^n}
\e\cases \f{7}{16}x^2-\f{p}{8}-\f{23p^2}{128x^2}\mod {p^3}&\t{if
$p=x^2+2y^2\e 1,3\mod 8$,}
\\-\f 3{64}R_2(p)-\f 3{32}p\mod{p^2}&\t{if $p\e 5,7\mod
8$}.\endcases$$
\endpro

\pro{Conjecture 3.34} Let $p$ be a prime with $p>3$. Then
$$\align &\Ls p3\sum_{n=0}^{p-1}n^2\f{V_n^{(6)}}{(-432)^n}
\\&\e\cases \f{5}{18}x^2-\f{5+18\sls p3}{36}p+\f{13p^2}{144x^2}\mod
{p^3}&\t{if $p=x^2+4y^2\e 1\mod 4$,}
\\-\f 1{24}R_1(p)+\f{1}{2}\Ls p3p\mod{p^2}&\t{if $p\e 3\mod
4$}\endcases\endalign$$ and
$$\align &\Ls p3\sum_{n=0}^{p-1}n^3\f{V_n^{(6)}}{(-432)^n}
\\&\e\cases \f{7}{12}x^2-\f{21-5\sls p3}{72}p-\f{19p^2}{96x^2}\mod
{p^3}&\t{if $p=x^2+4y^2\e 1\mod 4$,}
\\\f 1{16}R_1(p)-\f{5}{72}\Ls p3p\mod{p^2}&\t{if $p\e 3\mod
4$.}\endcases\endalign$$
\endpro
\par{\bf Remark 3.6} Let $p>3$ be a prime. For the congruences modulo $p^3$
for the sums
$$\align &\sum_{n=0}^{p-1}\f{V_n^{(3)}}{3^n},
\ \sum_{n=0}^{p-1}\f{V_n^{(3)}}{243^n},\
\sum_{n=0}^{p-1}\f{V_n^{(3)}}{(-27)^n},
\\&\sum_{n=0}^{p-1}V_n^{(4)},\
\sum_{n=0}^{p-1}\f{V_n^{(4)}}{4096^n}, \
\sum_{n=0}^{p-1}\f{V_n^{(4)}}{16^n}, \
 \sum_{n=0}^{p-1}\f{V_n^{(4)}}{256^n}, \
\\&\sum_{n=0}^{p-1}\f{V_n^{(4)}}{(-8)^n},\
\sum_{n=0}^{p-1}\f{V_n^{(4)}}{(-512)^n}, \
\sum_{n=0}^{p-1}\f{V_n^{(4)}}{(-64)^n}, \
\sum_{n=0}^{p-1}\f{V_n^{(6)}}{(-432)^n}
\endalign$$ see [S9].

\end{document}